\documentclass{article}

\usepackage{graphicx}
\usepackage{subfigure}

\newtheorem{thm}{Theorem}[section]
\newtheorem{cor}[thm]{Corollary}
\newtheorem{lem}[thm]{Lemma}
\newtheorem{de}[thm]{Definition}
\newcounter{bean}

\begin{document}

\title{Evenly Spaced Data Points and Radial Basis Functions}         
\author{Lin-Tian Luh\\Department of Mathematics, Providence University\\Shalu, Taichung, Taiwan\\Email:ltluh@pu.edu.tw\\phone:(04)26328001 ext. 15126 \\ fax:(04)26324653}          
\date{\today}          
\maketitle
{\bf Abstract}.This is a continuation of our study about shape parameter, based on an approach very different from that of \cite{Lu3} and \cite{Lu4}. Here we adopt an error bound of convergence order $O(\sqrt{d}\omega^{\frac{1}{d}})$ as $d\rightarrow 0$, where $0<\omega<1$ is a constant and $d$ denotes essentially the fill-distance. The constant $\omega$ is much smaller than the one appears in \cite{Lu3} and \cite{Lu4} where the error bound is $O(\omega^{\frac{1}{d}})$ only. Moreover, the constant $\omega$ here only mildly depends on the dimension $n$. It means that for high-dimensional problems the criteria of choosing the shape parameter presented in this paper are much better than those of \cite{Lu3} and \cite{Lu4}. The drawback is that the distribution of data points must be slightly controlled. 
\section{Introduction}       
We first review some basic material. The radial function we use to construct approximating functions is called shifted surface spline defined by
\begin{equation}
h(x):=(-1)^{m}(|x|^{2}+c^{2})^{\frac{\lambda}{2}}\log{(|x|^{2}+c^{2})^{\frac{1}{2}}},\ \lambda\in Z_{+},\ m=1+\frac{\lambda}{2},\ c>0,\ x\in R^{n},\ \lambda,n\ even,
\end{equation}
where $|x|$ is the Euclidean norm of $x$, log denotes the natural logarithm, and $\lambda,\ c$ are constants. The constant $c$ is called shape parameter which greatly influences the quality of the approximation. Unfortunately, its optimal choice is a big problem and has been regarded as a hard question, not only in mathematics, but also in engineering.

As is well known in the field of RBF(radial basis functions), for any scattered set of data points $(x_{1},f(x_{1})),\cdots, (x_{N},f(x_{N}))$, there is a unique function
\begin{equation}
s(x):=\sum_{j=1}^{N}c_{j}h(x-x_{j})+p(x)
\end{equation}
interpolating these data points, where $c_{1},\cdots,c_{N}$ are constants to be determined and $p(x)$ is a polynomial of degree $\leq m-1$. The only requirement for the data points is that $x_{1},\cdots,x_{N}$ should be polynomially nondegenerate.

The choice of $c$ severely influences the upper bound of $|f(x)-s(x)|$. To the author's knowledge, there are three kinds of error bound for shifted-surface-spline interpolation:algebraic type, exponential type, and improved exponential type. Among them the algebraic type shows nothing about the effect of $c$. The exponential type works well for this purpose, as can be seen in \cite{Lu3} and \cite{Lu4}. However the improved exponential type works better as will be seen in this article.
\subsection{Function Spaces}
We put restrictions on the approximated functions. 
\begin{de}
 For any $\sigma>0$, the class of band-limited functions $f$ in $L^{2}(R^{n})$ is defined by
$$B_{\sigma}:=\{ f\in L^{2}(R^{n}):\ \hat{f}(\xi)=0\ if\ |\xi|>\sigma\},$$
where $\hat{f}$ denotes the Fourier transform of $f$.
\end{de}
A larger function space is defined as follows.
\begin{de}
 For any $\sigma>0$,
$$E_{\sigma}:=\{ f\in L^{2}(R^{n}):\ \int|\hat{f}(\xi)|^{2}e^{\frac{|\xi|^{2}}{\sigma}}d\xi<\infty\},$$
where $\hat{f}$ denotes the Fourier transform of $f$. For each $f\in E_{\sigma}$, its norm is 
$$\|f\|_{E_{\sigma}}:=\left\{\int|\hat{f}(\xi)|^{2}e^{\frac{|\xi|^{2}}{\sigma}}d\xi\right\}^{\frac{1}{2}}.$$
\end{de}
Although we only deal with functions from $B_{\sigma}$ or $E_{\sigma}$ in this article, another function space should be mentioned. It's denoted by ${\cal C}_{h,m}$ and is the so-called native space induced by $h$. We omit its complicated definition and characterization. For these, we refer the reader to \cite{Lu1,Lu2,MN1,MN2,We}. What we need here is just $B_{\sigma}\subseteq E_{\sigma}\subseteq {\cal C}_{h,m}$. The proof of $E_{\sigma}\subseteq {\cal C}_{h,m}$ can be found in \cite{Lu4}. Moreover, each function $f\in {\cal C}_{h,m}$ has a semi-norm denoted by $\|f\|_{h}$.
\subsection{Distribution of Data Points}
The distribution of data points plays a crucial role in our approach. We first review a basic definition of \cite{Lu4}.
\begin{de}
Let E be an n-dimensional simplex in $R^{n}$ with vertices $v_{1},\cdots, v_{n+1}$. For any point $x\in E$, its barycentric coordinates are the numbers $\lambda_{1},\cdots, \lambda_{n+1}$ satisfying
$$x=\sum_{i=1}^{n+1}\lambda_{i}v_{i},\ \sum_{i=1}^{n+1}\lambda_{i}=1,\ \lambda_{i}\geq 0\ for\ all\ i.$$
\end{de}
The definition of simplex can be found in \cite{Fl}.
\begin{de}
For any n-dimensional simplex, the evenly spaced points of degree k are the points whose barycentric coordinates are of the form $$(k_{1}/k,k_{2}/k,\cdots,k_{n+1}/k),\ k_{i}\ nonnegative\ integers\ and\ k_{1}+\cdots,k_{n+1}=k.$$  
\end{de}
As pointed out in \cite{Lu5}, the number of such points is equal to the dimension of $P_{k}^{n}$, the space of n-dimensional polynomials of degree less than or equal to $k$. In this paper we use $N$ to denote this number. Thus $N=dimP_{k}^{n}$.

In our approach interpolation occurs in a simplex and the centers(interpolation points) are evenly spaced points of that simplex. Note that the shape of the simplex is very flexible and hence this requirement is not very restrictive. We shall see in the next section that for this kind of interpolation there is an error bound which is much better than the case of purely scattered data points, making the criteria of choosing $c$ much more meaningful.
\section{Improved Exponential-type Error Bound}
The function $h(x)$ in (1) induces a few basic ingredients of the error bound. First, its Fourier transform is
\begin{eqnarray}
\hat{h}(\xi)=l(\lambda,n)|\xi|^{-\lambda-n}\tilde{{\cal K}}_{\frac{n+\lambda}{2}}(c|\xi|)
\end{eqnarray}
where $l(\lambda,n)$ is a constant depending on $\lambda,\ n$\cite{Lu5,MN2}, and $\tilde{{\cal K}}_{\nu}(t)=t^{\nu}{\cal K}_{\nu}(t),\ {\cal K}_{\nu}(t)$ being the modified Bessel function of the second kind \cite{AS}. Second, each $h(x)$ corresponds to two constants $\rho$ and $\Delta_{0}$ defined as follows.
\begin{de}
 Let $h(x)$ be as in (1). The constants $\rho$ and $\Delta_{0}$ are defined as follows.
\begin{list}
  {(\alph{bean})}{\usecounter{bean} \setlength{\rightmargin}{\leftmargin}}
  \item Suppose $n-\lambda >3$. Let $s=\lceil \frac{n-\lambda -3}{2}\rceil $. Then 
      $$\rho=1+\frac{s}{2m+3} \ and \ \Delta_{0}=\frac{(2m+2+s)(2m+1+s)\cdots (2m+3)}{\rho^{2m+2}}.$$     
  \item Suppose $n-\lambda \leq 1$. Let $s=-\lceil \frac{n-\lambda-3}{2}\rceil$. Then
   $$\rho=1\ and\ \Delta_{0}=\frac{1}{(2m+2)(2m+1)\cdots (2m-s+3)}.$$
  \item Suppose $1<n-\lambda\leq 3$. Then
 $$\rho =1\ and \ \Delta_{0}=1.$$  
\end{list}
\end{de}
The following core theorem provides the theoretical ground of our criteria of choosing $c$. We cite it directly from \cite{Lu5}. with a slight modification.
\begin{thm}
 Let $h$ be as in (1). For any positive number $b_{0}$, there exist positive constants $\delta_{0},\ c_{1},\ C,\ \omega,\ 0<\omega<1$, completely determined by $h$ and $b_{0}$, such that for any n-dimensional simlex $Q_{0}$ of diameter $b_{0}$, any $f\in {\cal C}_{h,m}$, and any $0<\delta\leq \delta_{0}$, there is a number $r$ satisfying the property that $\frac{1}{3C}\leq r\leq b_{0}$ and for any n-dimensional simplex Q of diameter r, $Q\subseteq Q_{0}$, there is an interpolating function $s(\cdot )$ as defined in (2) such that
\begin{equation}
 |f(x)-s(x)|\leq c_{1}\sqrt{\delta}(\omega)^{\frac{1}{\delta}}\|f\|_{h}
\end{equation}
for all $x$ in Q, where $C$ is defined by
$$C:=\max\left\{ 8\rho',\ \frac{2}{3b_{0}}\right\},\ \rho'=\frac{\rho}{c}$$
where $\rho$ and $c$ appeared in Definition2.1 and (1) respectively. The function $s(\cdot )$ interpolates $f$ at $x_{1},\cdots ,x_{N}$ which are evenly spaced points of degree $k-1$ on $Q$ as defined in Definition1.4, with $k=\frac{r}{\delta}$. Here $\|f\|_{h}$ is the h-norm of $f$ in the native space.

The numbers $\delta_{0},\ c_{1}$ and $\omega$ are given by $\delta_{0}:=\frac{1}{3C(m+1)}$ where m appeared in (1), $$c_{1}:=\sqrt{l(\lambda,n)}(2\pi)^{\frac{1}{4}}\sqrt{n\alpha_{n}}c^{\frac{\lambda}{2}}\sqrt{\Delta_{0}}\sqrt{3C}\sqrt{(16\pi)^{-1}}$$ where $\lambda$ is as in (1), $l(\lambda,n)$ appeared in (3), $\alpha_{n}$ is the volume of the unit ball in $R^{n}$, and $\Delta_{0}$ was defined in Definition2.1, and $\omega:=\left(\frac{2}{3}\right)^{\frac{1}{3C}}$. 
\end{thm}
{\bf Remark}: This seemingly complicated theorem is in fact not so difficult to understand. The number $\delta$ is in spirit the well-known fill-distance.\\
\\
Now, it's easily seen in (4) both $c_{1}$ and $\omega$ depend on the shape parameter $c$. So does $\|f\|_{h}$. In order to make (4) useful for choosing $c$, we still have to convert $\|f\|_{h}$ into a transparent expression of $c$. We need two lemmas which we cite from \cite{Lu3} and \cite{Lu4}, respectively.
\begin{lem}
 For any $\sigma>0$, $f\in B_{\sigma}$ implies $f\in {\cal C}_{h,m}$ and 
$$\|f\|_{h}\leq C_{0}(m,n)\cdot \left(\frac{2}{\pi}\right)^{\frac{1}{4}}\cdot \sigma^{\frac{1+n+\lambda}{4}}\cdot c^{\frac{1-n-\lambda}{4}}\cdot e^{\frac{c\sigma}{2}}\cdot\|f\|_{L^{2}}$$
, where
$$C_{0}(m,n)=\frac{(2\pi)^{-n}\sqrt{m!}}{\sqrt{l(\lambda,n)}}.$$
\end{lem} 
\begin{lem}
For any $\sigma>0,\ f\in E_{\sigma}$ implies $f\in {\cal C}_{h,m}$ and
$$\|f\|_{h}\leq a_{0}(\lambda,n)c^{\frac{1-n-\lambda}{4}}\sup_{\xi\in R^{n}}\left\{|\xi|^{\frac{1+n+\lambda}{4}}e^{\frac{c|\xi|}{2}-\frac{|\xi|^{2}}{2\sigma}}\right\}\|f\|_{E_{\sigma}}$$
, where $a_{0}(\lambda,n)=\sqrt{\frac{m!}{l(\lambda,n)}}2^{\frac{1}{4}-n}\pi^{-n-\frac{1}{4}}$. 
\end{lem} 
\begin{cor}
If $f\in B_{\sigma}$, (4) can be converted into
\begin{equation}
 |f(x)-s(x)|\leq C_{B}c^{\frac{1+\lambda-n}{4}}\sqrt{C}e^{\frac{c\sigma}{2}}\left(\frac{2}{3}\right)^{\frac{1}{3C\delta}}\sqrt{\delta}\|f\|_{L^{2}}
\end{equation}
where $C_{B}:=\sigma^{\frac{1+n+\lambda}{4}}(2\pi)^{-n}\sqrt{6n\alpha_{n}\Delta_{0}(16\pi)^{-1}m!}$.
\end{cor}
\begin{cor}
 If $f\in E_{\sigma}$, (4) can be converted into
\begin{equation}
 |f(x)-s(x)|\leq C_{E}c^{\frac{1+\lambda-n}{4}}\sqrt{C}\sup_{\xi\in R^{n}}\left\{|\xi|^{\frac{1+n+\lambda}{4}}e^{\frac{c|\xi|}{2}-\frac{|\xi|^{2}}{2\sigma}}\right\}\left(\frac{2}{3}\right)^{\frac{1}{3C\delta}}\sqrt{\delta}\|f\|_{E_{\sigma}}
\end{equation}
where $C_{E}:=(2\pi)^{-n}\sqrt{6n\alpha_{n}\Delta_{0}(16\pi)^{-1}m!}$.
\end{cor}
\section{Criteria of Choosing $c$}
Note that in the right-hand side of both (5) and (6), after every thing independent of $c$ is fixed, there is a function of $c$ which may be used to choose the optimal $c$. However, in Theorem2.2 there is still an object dependent of $c$. That is $\delta_{0}$, the upper bound of $\delta$. For a fixed $\delta>0$, decreasing $c$ may decrease $\delta_{0}$ and the requirement $\delta\leq \delta_{0}$ may be violated. For any fixed $b_{0}>0$, the minimal acceptable $c$ is $c=c_{0}:=24\rho(m+1)\delta$. Also, $C=\frac{2}{3b_{0}}$ if and only if $c\geq c_{1}:=12\rho b_{0}$. Here $c_{1}$ is different from the $c_{1}$ in (4).

There is a logical problem about $c,\ b_{0},\ \delta$ and $\delta_{0}$. According to Theorem2.2, $c$ appears first, then $b_{0}$, and then $\delta_{0}$ and $\delta$. There is a trick to resolve this logical problem. For any $b_{0}>0$, we first let $C=\frac{2}{3b_{0}}$ and $\delta_{0}:=\frac{1}{3C(m+1)}=\frac{b_{0}}{2(m+1)}$ temporarily. Let $\delta<\delta_{0}$ be fixed. Then $\delta\leq \delta_{0}$ will always be satisfied if $c\in[c_{0},\infty)$. After the optimal $c$ is obtained, we let $C$ and $\delta_{0}$ be as in Theorem2.2. The consequence is that we can only choose an optimal $c$ from $[c_{0},\infty)$. This is a drawback of our theory. However, since $\delta$ can be theoretically arbitrarily small, $c_{0}$ can be very close to zero, theoretically.

Now, on the right side of either (5) or (6), there is a function of $c$. Let's call it an MN function and denote it by $MN(c)$. Thus for $f\in B_{\sigma}$,
\begin{equation}
 MN(c)=c^{\frac{1+\lambda-n}{4}}\sqrt{C}e^{\frac{c\sigma}{2}}\left(\frac{2}{3}\right)^{\frac{1}{3C\delta}}
\end{equation}
, and for $f\in E_{\sigma}$,
\begin{equation}
 MN(c)=c^{\frac{1+\lambda-n}{4}}\sqrt{C}\sup_{\xi \in R^{n}}\left\{ |\xi|^{\frac{1+n+\lambda}{4}}e^{\frac{c|\xi|}{2}-\frac{|\xi|^{2}}{2\sigma}}\right\}\left(\frac{2}{3}\right)^{\frac{1}{3C\delta}}
\end{equation}
. Our goal is to find $c\in [c_{0},\infty)$ which minimizes $MN(c)$.

In the preceding discussion the domain size $b_{0}$ was fixed. For $b_{0}$ not fixed, the way of choosing $c$ will be different. We deal with them separately.
\subsection{$b_{0}$ fixed}
Note that 
$$C=\left\{ \begin{array}{ll}
              \frac{8\rho}{c} & \mbox{if $c_{0}\leq c\leq c_{1}$},\\ 
             \\ 
              \frac{2}{3b_{0}} & \mbox{if $c_{1}\leq c<\infty$}
            \end{array} \right. $$
, where $c_{0}$ and $c_{1}$ were defined in the beginning of section3. Thus (7) and (8) can be refined as
\begin{equation}
MN(c)=\left\{ \begin{array}{ll}
               \sqrt{8\rho}c^{\frac{\lambda-n-1}{4}}e^{c\left[\frac{\sigma}{2}+\frac{\ln{\frac{2}{3}}}{24\rho \delta}\right]} & \mbox{if $c_{0}\leq c\leq c_{1}$}\\
                                                    \\  
               \sqrt{\frac{2}{3b_{0}}}c^{\frac{1+\lambda-n}{4}}e^{\frac{c\sigma}{2}}\left(\frac{2}{3}\right)^{\frac{b_{0}}{2\delta}} & \mbox{if $c_{1}\leq c<\infty$}
              \end{array} \right.
\end{equation}
for $f\in B_{\sigma}$, and
\begin{equation}
  MN(c)=\left\{ \begin{array}{ll}
                  \sqrt{8\rho}c^{\frac{\lambda-n-1}{4}}\sup_{\xi\in R^{n}}\left\{ |\xi|^{\frac{1+n+\lambda}{4}}e^{\frac{c|\xi|}{2}-\frac{|\xi|^{2}}{2\sigma}}\right\}\left(\frac{2}{3}\right)^{\frac{c}{24\rho\delta}} & \mbox{if $c_{0}\leq c\leq c_{1}$} \\
                                       \\
                  \sqrt{\frac{2}{3b_{0}}}c^{\frac{1+\lambda-n}{4}}\sup_{\xi\in R^{n}}\left\{|\xi|^{\frac{1+n+\lambda}{4}}e^{\frac{c|\xi|}{2}-\frac{|\xi|^{2}}{2\sigma}}\right\}\left(\frac{2}{3}\right)^{\frac{b_{0}}{2\delta}} & \mbox{if $c_{1}\leq c<\infty$}
                \end{array} \right.
\end{equation}
for $f\in E_{\sigma}$.
\subsubsection{$f\in B_{\sigma}$}
For $f\in B_{\sigma}$, we have the following cases, where $k:=\frac{\sigma}{2}+\frac{\ln{\frac{2}{3}}}{24\rho\delta}$.\\
\\
{\bf Case1}. \fbox{$\lambda-n-1\geq 0$} and \fbox{$k\geq 0$} For any $b_{0}>0$ in Theorem2.2 and positive $\delta<\frac{b_{0}}{2(m+1)}$, if $\lambda-n-1\geq 0$ and $k\geq 0$, the optimal choice of $c$ for $c\in [c_{0},\infty)$ is to let $c=c_{0}:=24\rho(m+1)\delta$.\\
\\
{\bf Reason}: In this case $MN(c)$ in (9) is increasing on $[c_{0},\infty)$. \\
\\
{\bf Numerical Examples}:
\begin{figure}[h]
\centering
\mbox{
\subfigure[a smaller domain]{\includegraphics[scale=0.9]{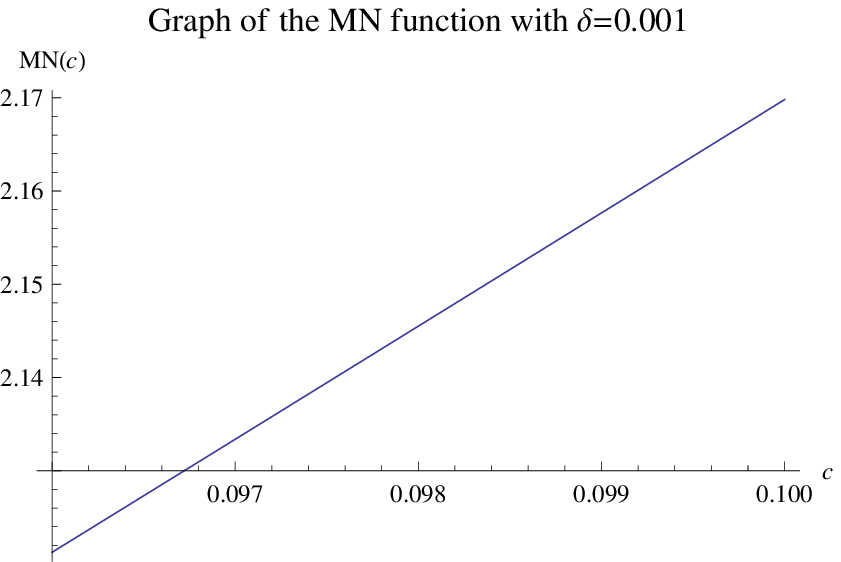}}
\subfigure[a larger domain]{\includegraphics[scale=0.9]{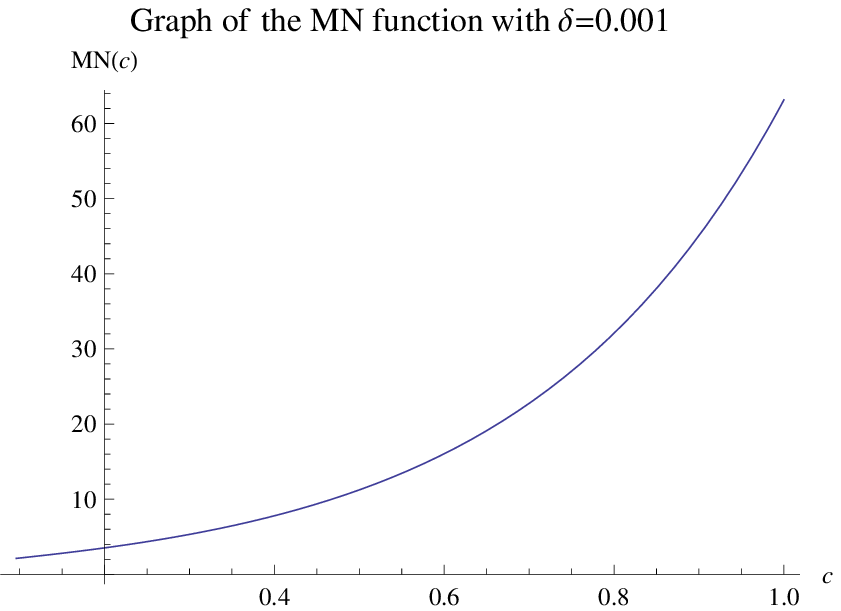}}}
\caption{Here $n=2,\lambda=4,\sigma=40$ and $b_{0}=10$.}
\end{figure}
\\
{\bf Case2}. \fbox{$\lambda-n-1\geq 0$} and \fbox{$k<0$} For any $b_{0}>0$ in Theorem2.2 and positive $\delta<\frac{b_{0}}{2(m+1)}$, if $\lambda-n-1\geq 0$ and $k<0$, the optimal choice of $c$ for $c\in [c_{0}, \infty)$ is $c^{*}\in [c_{0},c_{1}]$ which minimizes $MN(c)$ of (9) on $[c_{0},c_{1}]$.\\
\\
{\bf Reason}: In this case $MN(c)$ in (9) is increasing on $[c_{1},\infty)$.\\
\\
{\bf Numerical Examples}:

\begin{figure}[t]
\centering
\mbox{
\subfigure[a smaller domain]{\includegraphics[scale=0.9]{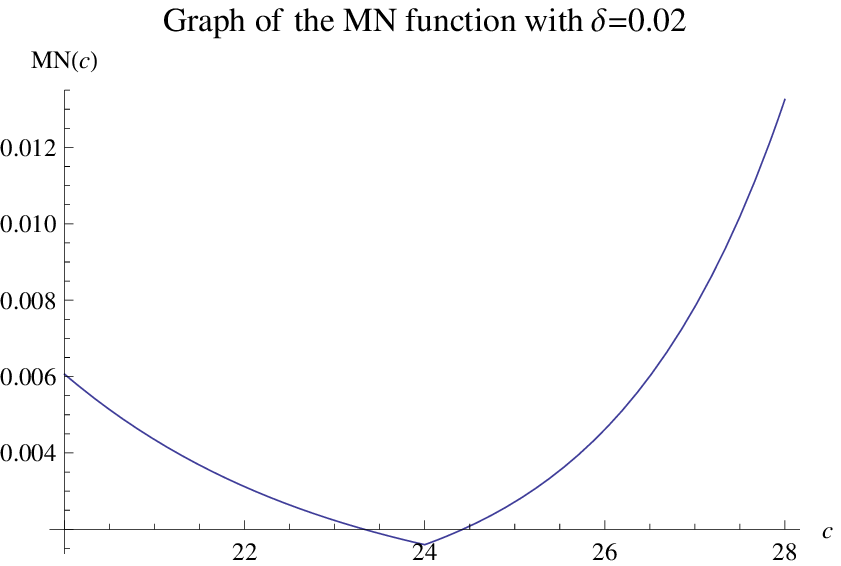}}
\subfigure[a larger domain]{\includegraphics[scale=0.9]{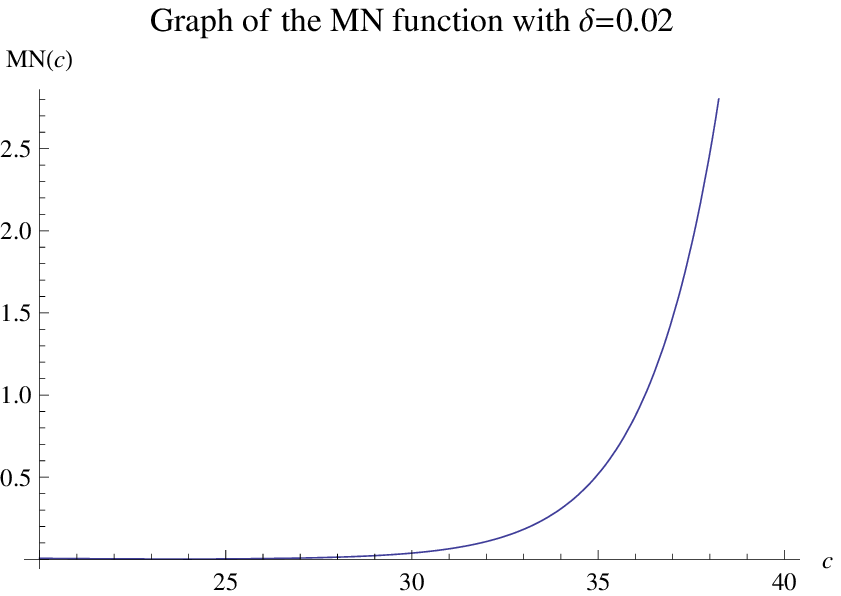}}}
\caption{Here $n=2,\lambda=4,\sigma=1$ and $b_{0}=2$.}

\mbox{
\subfigure[a smaller domain]{\includegraphics[scale=0.9]{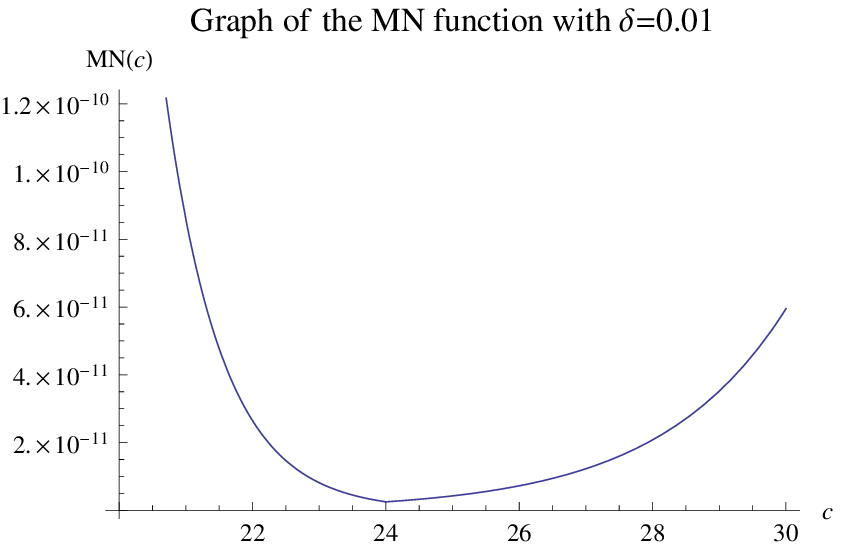}}
\subfigure[a larger domain]{\includegraphics[scale=0.9]{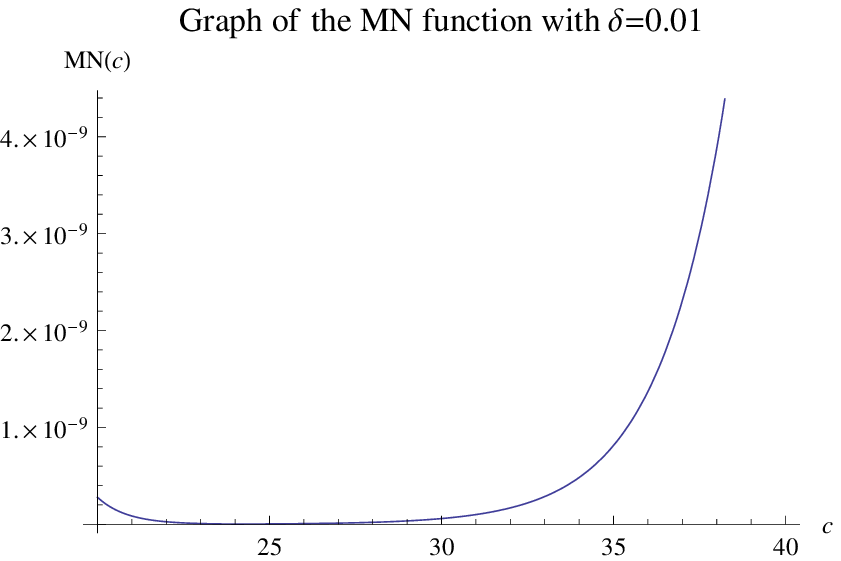}}}
\caption{Here $n=2,\lambda=4,\sigma=1$ and $b_{0}=2$.}

\mbox{
\subfigure[a smaller domain]{\includegraphics[scale=0.9]{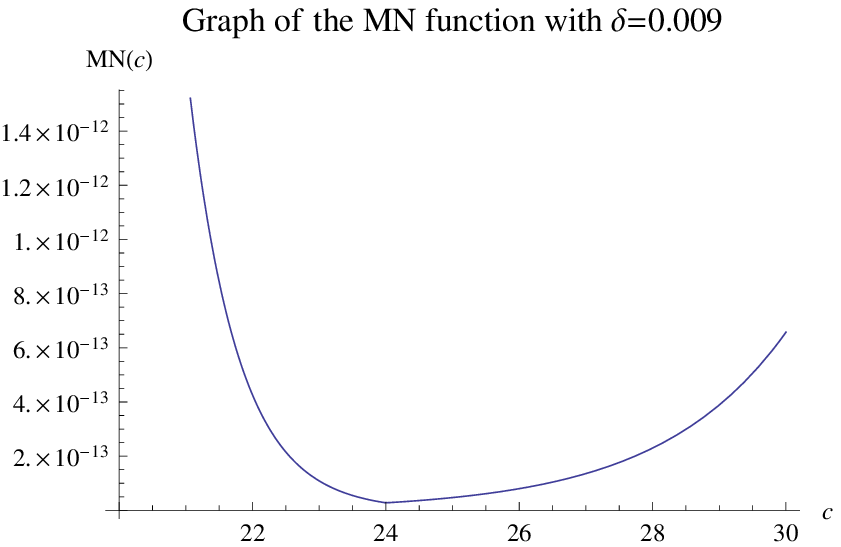}}
\subfigure[a larger domain]{\includegraphics[scale=0.9]{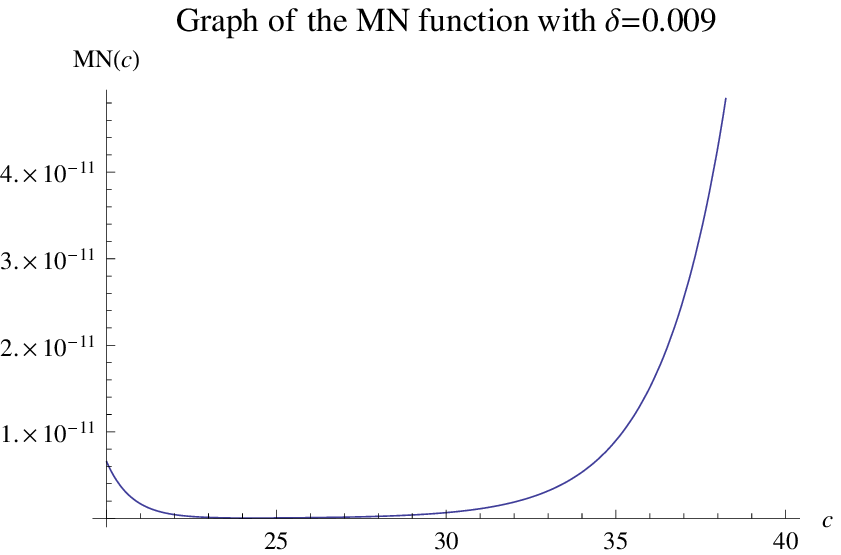}}}
\caption{Here $n=2,\lambda=4,\sigma=1$ and $b_{0}=2$.}

\end{figure}

\clearpage

\begin{figure}[t]
\centering
\mbox{
\subfigure[a smaller domain]{\includegraphics[scale=0.9]{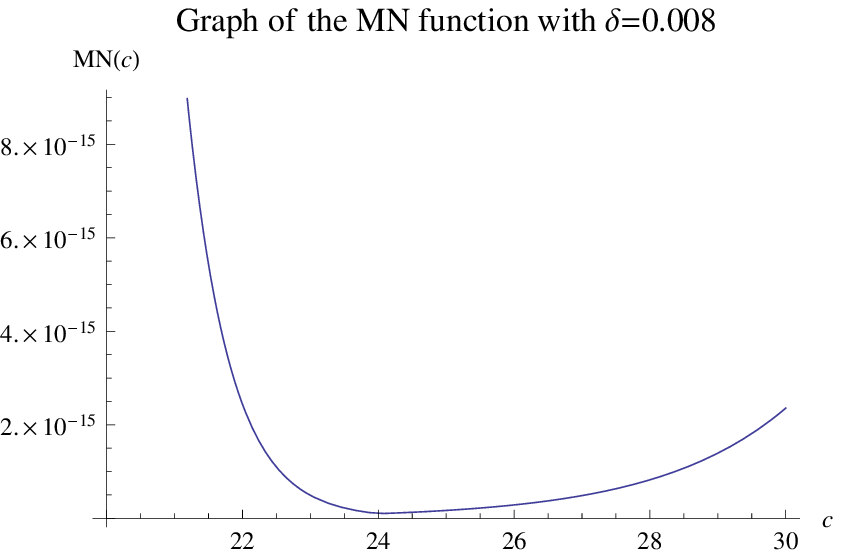}}
\subfigure[a larger domain]{\includegraphics[scale=0.9]{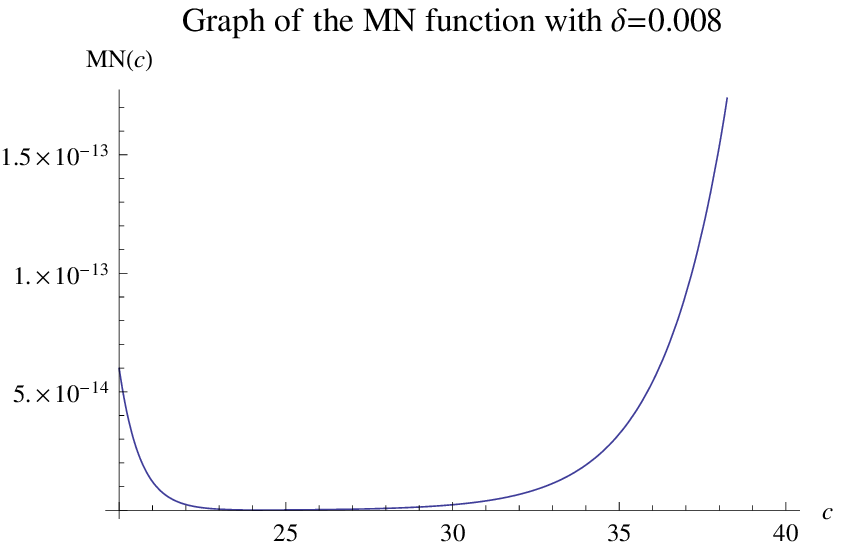}}}
\caption{Here $n=2,\lambda=4,\sigma=1$ and $b_{0}=2$.}

\mbox{
\subfigure[a smaller domain]{\includegraphics[scale=0.9]{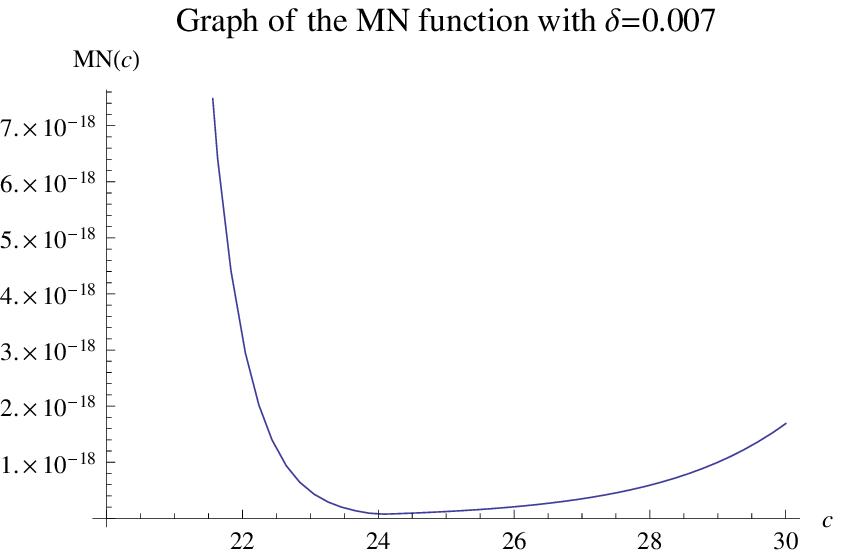}}
\subfigure[a larger domain]{\includegraphics[scale=0.9]{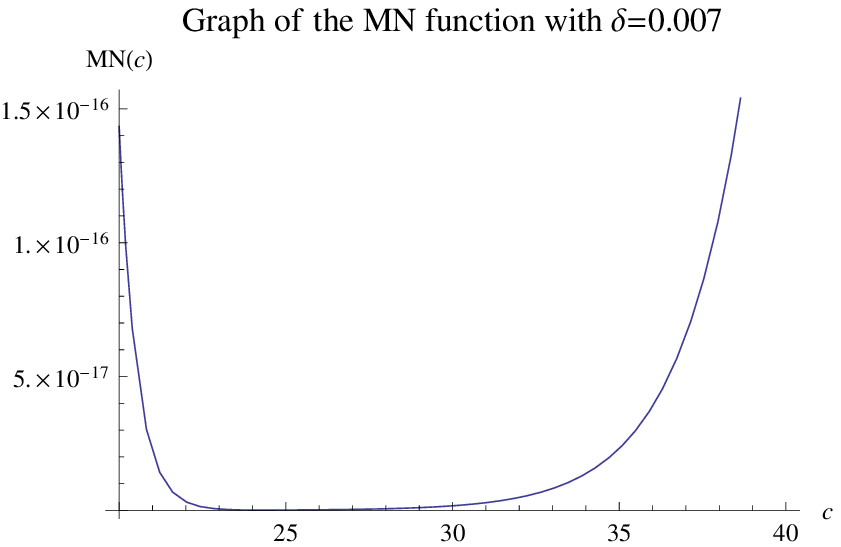}}}
\caption{Here $n=2,\lambda=4,\sigma=1$ and $b_{0}=2$.}

\end{figure}
All these figures provide a very small part of the entire curve only. In fact the curve increases or decreases very rapidly whenever $c$ is far from its optimal value.\\
\\ 
{\bf Case3}. \fbox{$\lambda-n-1<0$} and \fbox{$k<0$} For any $b_{0}>0$ in Theorem2.2 and positive $\delta<\frac{b_{0}}{2(m+1)}$, if $\lambda-n-1<0$ and $k<0$, the optimal choice of $c\in [c_{0},\infty)$ is the value $c^{*}\in [c_{1},\infty)$ which minimizes $MN(c)$ in (9) on $[c_{1},\infty)$.\\
\\
{\bf Reason}: In this case $MN(c)$ in (9) decreases on $[c_{0},c_{1}]$.\\
\\
{\bf Remark}: In Case3 if $1+\lambda-n\geq 0,\ MN(c)$ will be increasing on $[c_{1},\infty)$ and $c^{*}=c_{1}$.\\
\\
{\bf Numerical Examples}:
\begin{figure}[t]
\centering
\mbox{
\subfigure[a smaller domain]{\includegraphics[scale=0.9]{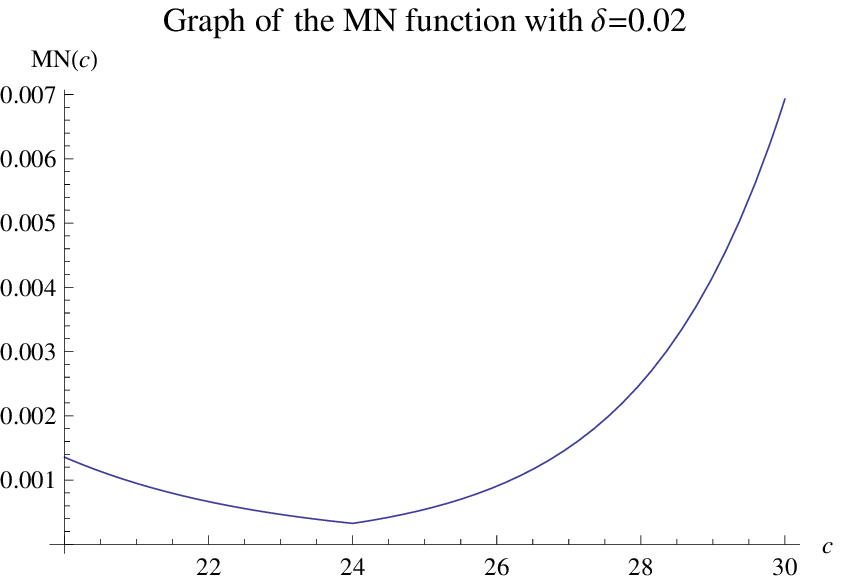}}
\subfigure[a larger domain]{\includegraphics[scale=0.9]{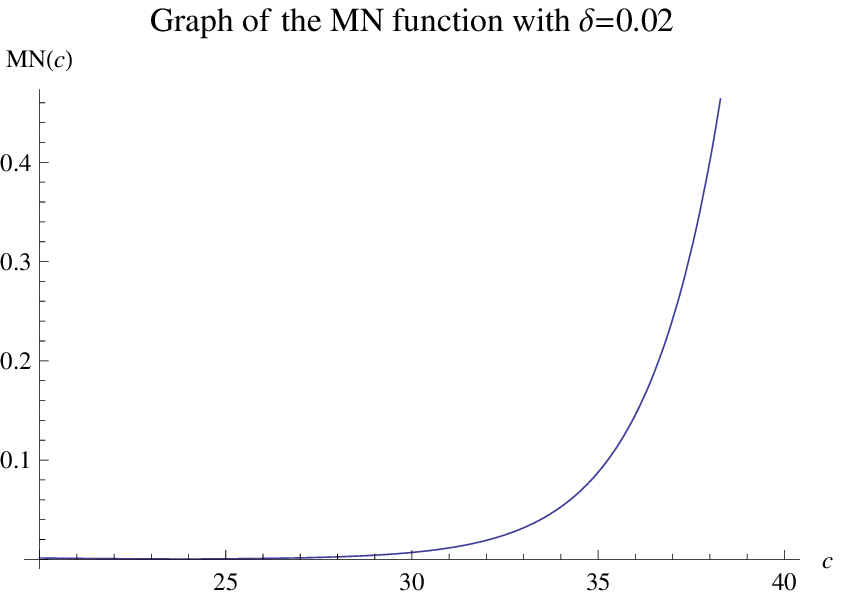}}}
\caption{Here $n=2,\lambda=2,\sigma=1$ and $b_{0}=2$.}

\mbox{
\subfigure[a smaller domain]{\includegraphics[scale=0.9]{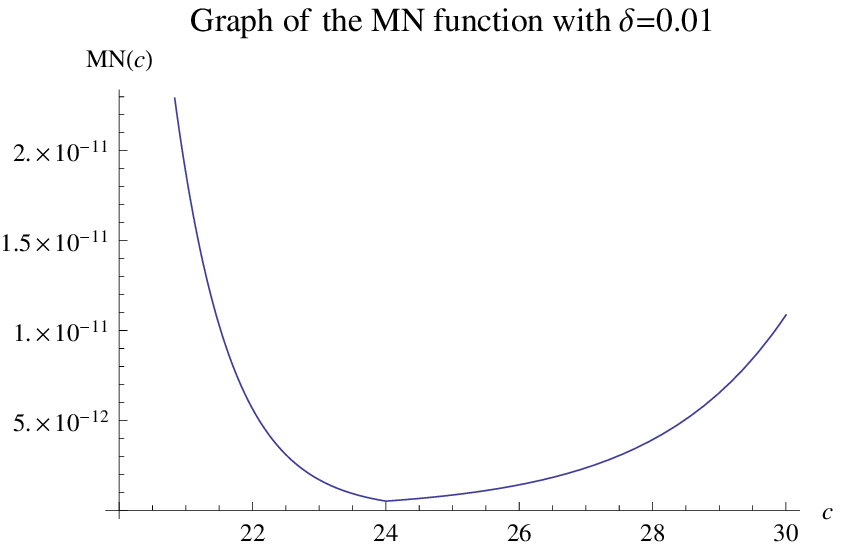}}
\subfigure[a larger domain]{\includegraphics[scale=0.9]{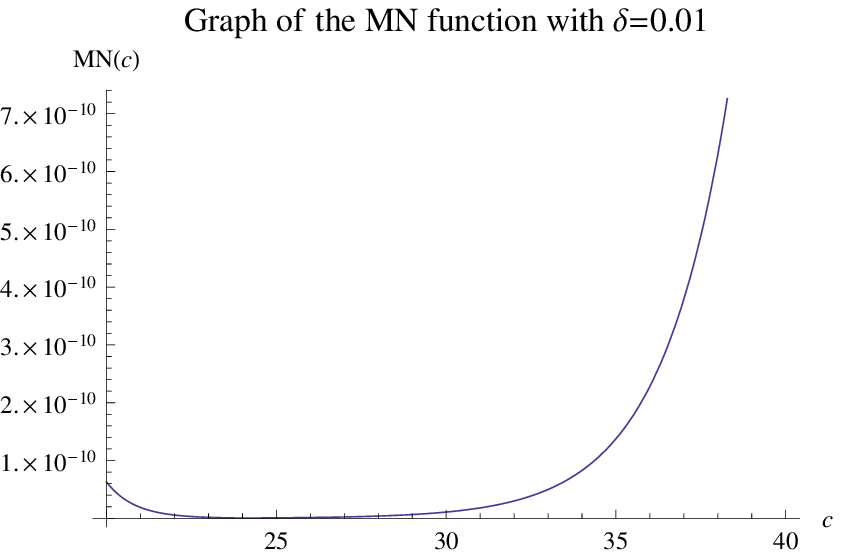}}}
\caption{Here $n=2,\lambda=2,\sigma=1$ and $b_{0}=2$.}

\mbox{
\subfigure[a smaller domain]{\includegraphics[scale=0.9]{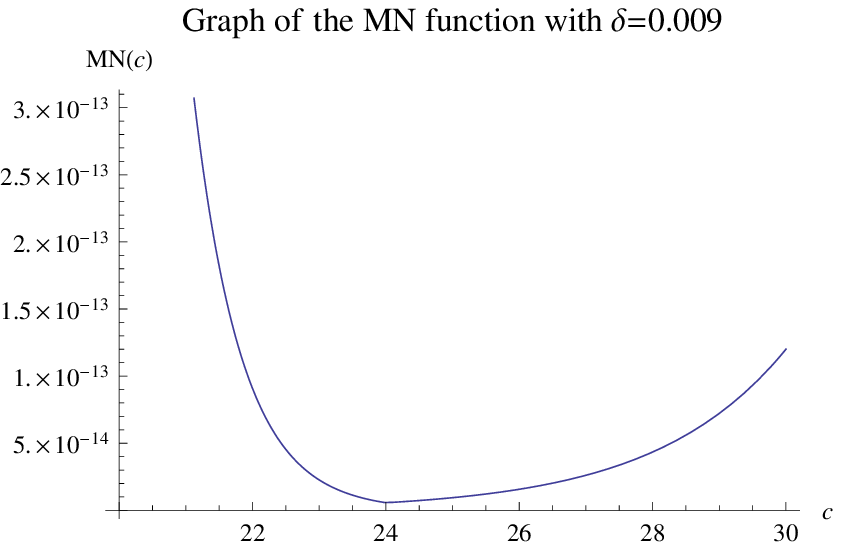}}
\subfigure[a larger domain]{\includegraphics[scale=0.9]{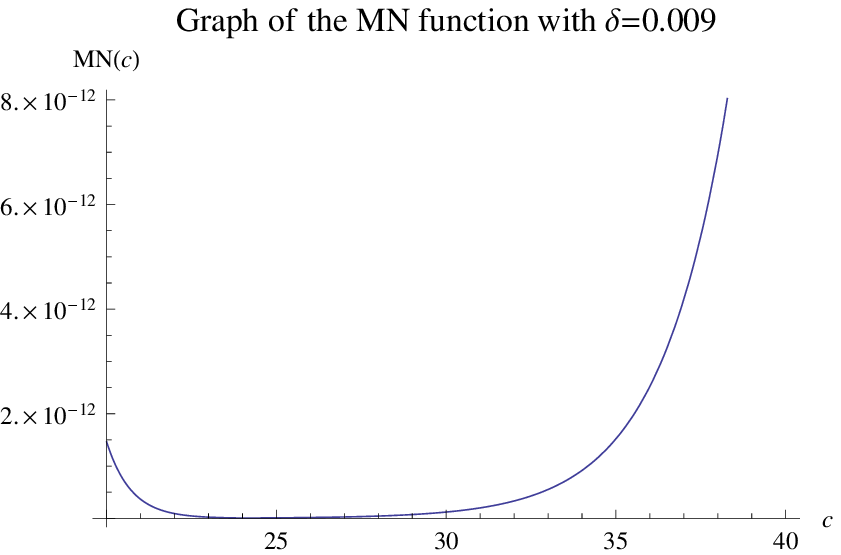}}}
\caption{Here $n=2,\lambda=2,\sigma=1$ and $b_{0}=2$.}
\end{figure}

\clearpage

\begin{figure}[t]
\centering
\mbox{
\subfigure[a smaller domain]{\includegraphics[scale=0.9]{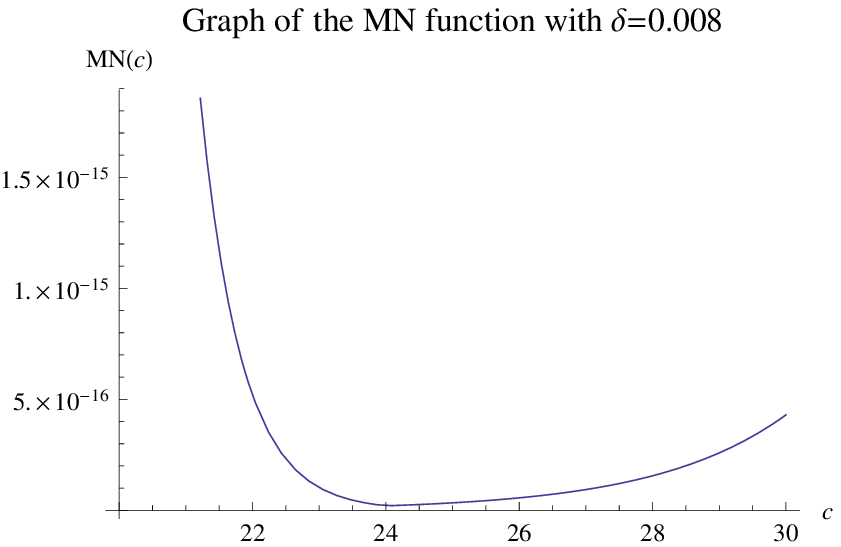}}
\subfigure[a larger domain]{\includegraphics[scale=0.9]{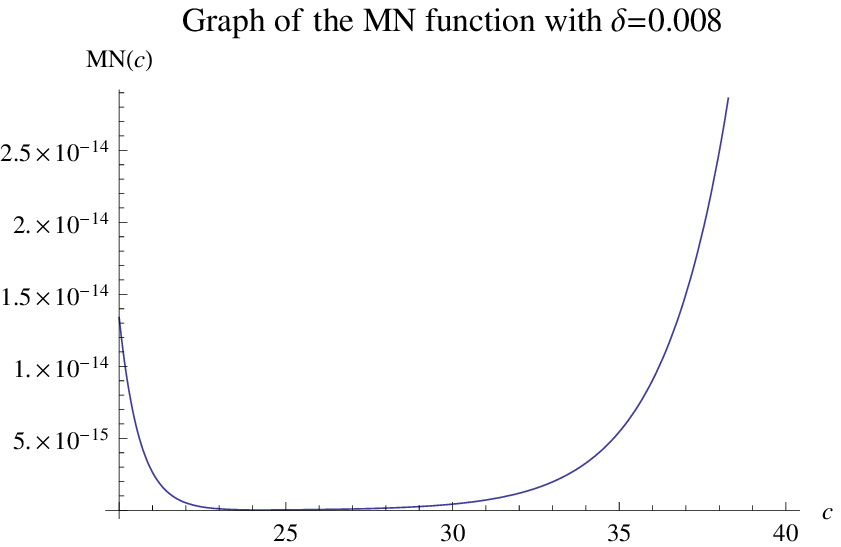}}}
\caption{Here $n=2,\lambda=2,\sigma=1$ and $b_{0}=2$.}

\mbox{
\subfigure[a smaller domain]{\includegraphics[scale=0.9]{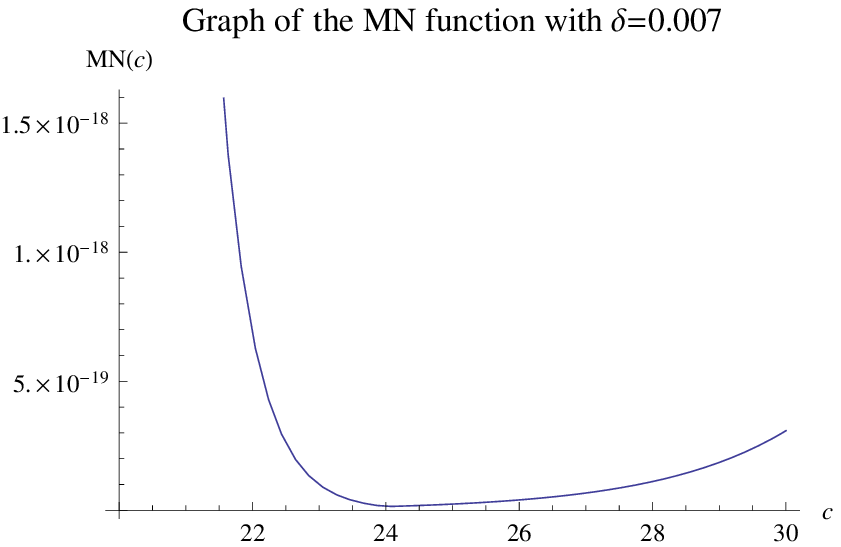}}
\subfigure[a larger domain]{\includegraphics[scale=0.9]{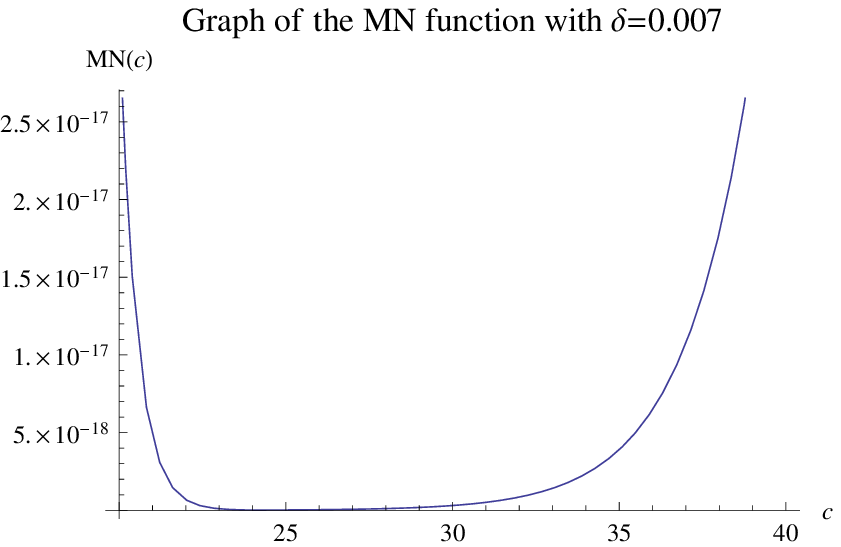}}}
\caption{Here $n=2,\lambda=2,\sigma=1$ and $b_{0}=2$.}
\end{figure}
The following case is quite different from the preceding three cases. For fixed $\sigma>0$, the number $\delta$ cannot be arbitrarily small due to the restriction $k>0$. However the optimal choice of $c$ depends on the domain size $b_{0}$.\\
\\
{\bf Case4}. \fbox{$\lambda -n-1<0$} and \fbox{$k\geq 0$} For any $b_{0}>0$ in Theorem2.2 and positive $\delta<\frac{b_{0}}{2(m+1)}$, if $\lambda-n-1<0$ and $k\geq 0$, the optimal choice of $c\in[c_{0},\infty)$ is either $c^{*}\in [c_{0},c_{1}]$ or $c^{**}\in [c_{1},\infty)$, depending on $MN(c^{*})\leq MN(c^{**})$ or $MN(c^{**})\leq MN(c^{*})$, where $c^{*}$ and $c^{**}$ minimize $MN(c)$ in (9) on $[c_{0},c_{1}]$ and $[c_{1},\infty)$, respectively.\\
\\
{\bf Reason}: In this case $MN(c)$ in (9) may not be monotonic on both $[c_{0},c_{1}]$ and $[c_{1},\infty)$.\\
\\
{\bf Remark}: In Case4 if $1+\lambda-n\geq 0,\ MN(c)$ will be increasing on $[c_{1},\infty)$ and $c^{**}=c_{1}$.\\
\\
{\bf Numerical Examples}:
\begin{figure}[t]
\centering
\mbox{
\subfigure[a smaller domain]{\includegraphics[scale=0.9]{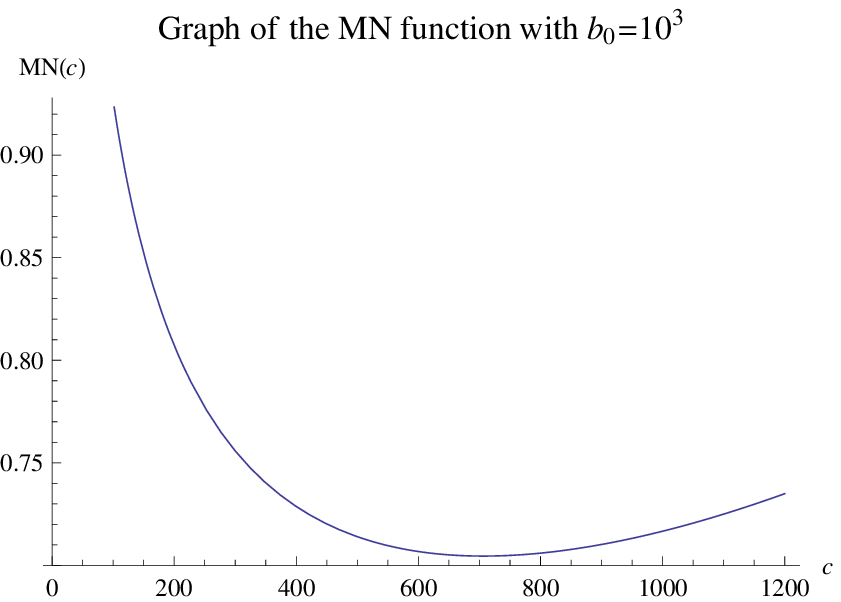}}
\subfigure[a larger domain]{\includegraphics[scale=0.9]{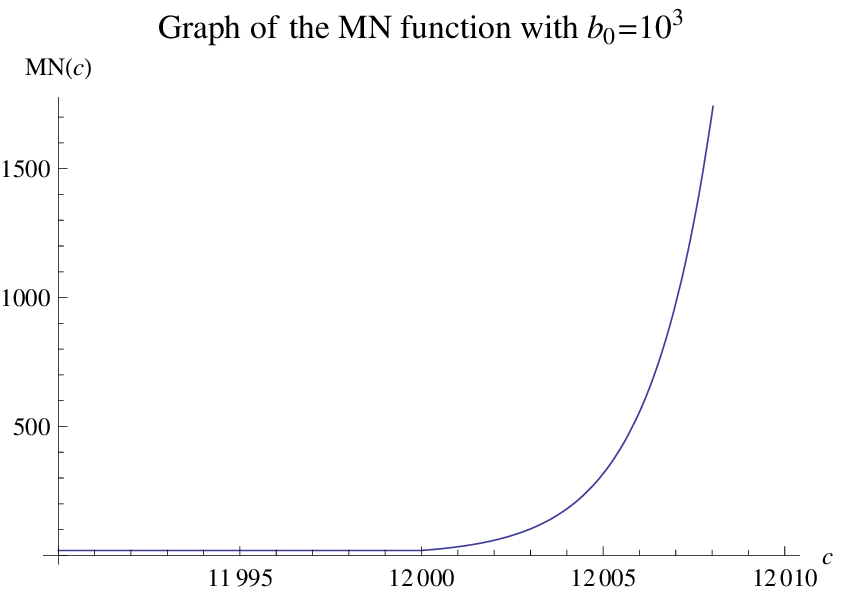}}}
\caption{Here $n=2,\lambda=2,\sigma=1.127$ and $\delta=0.03$.}

\mbox{
\subfigure[a smaller domain]{\includegraphics[scale=0.9]{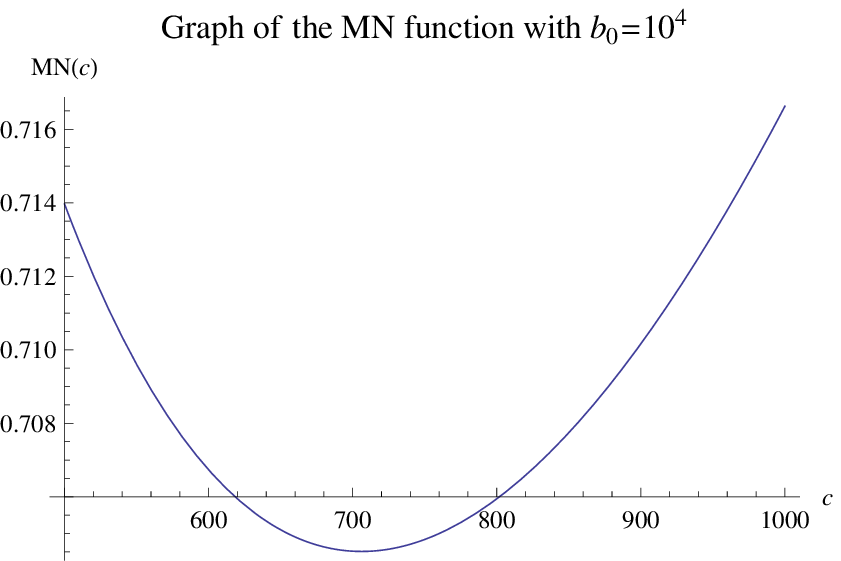}}
\subfigure[a larger domain]{\includegraphics[scale=0.9]{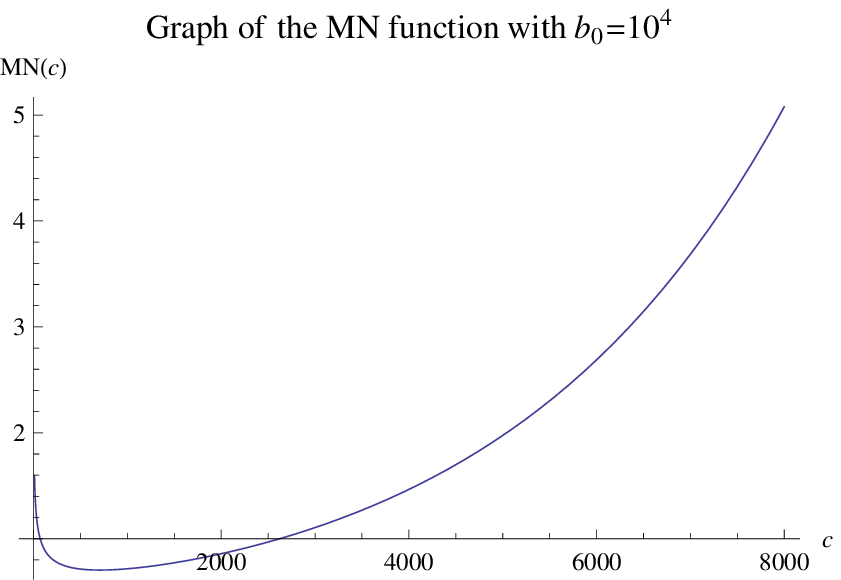}}}
\caption{Here $n=2,\lambda=2,\sigma=1.127$ and $\delta=0.03$.}

\mbox{
\subfigure[a smaller domain]{\includegraphics[scale=0.9]{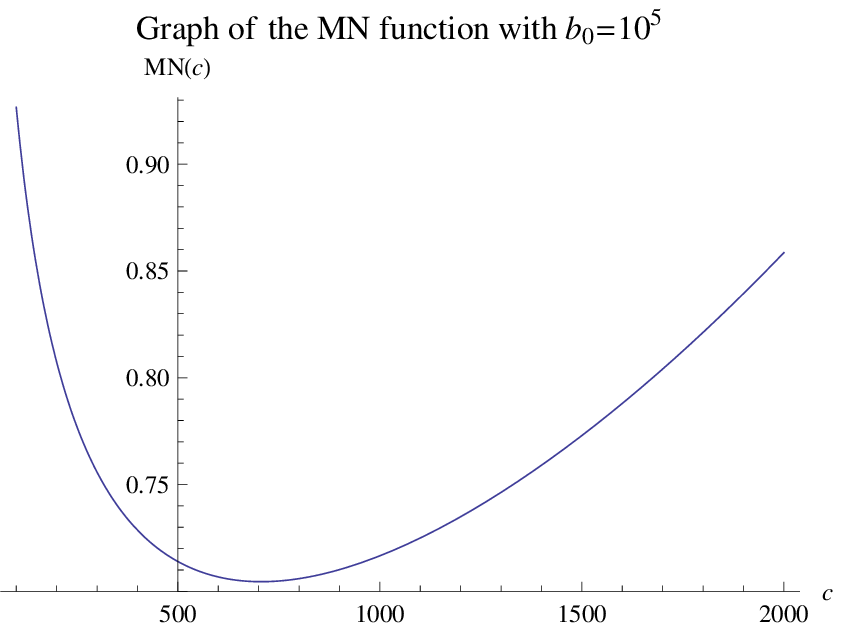}}
\subfigure[a larger domain]{\includegraphics[scale=0.9]{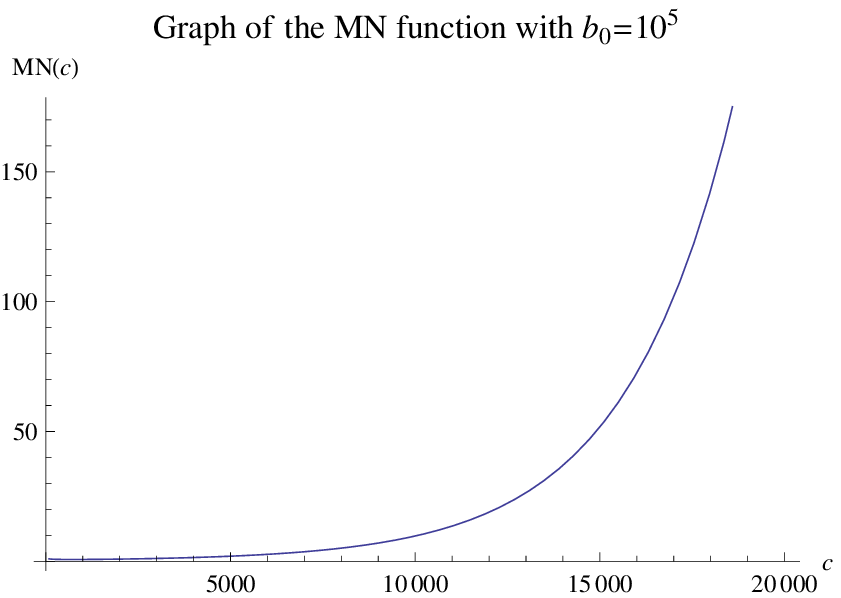}}}
\caption{Here $n=2,\lambda=2,\sigma=1.127$ and $\delta=0.03$.}
\end{figure}

\clearpage

\begin{figure}[t]
\centering
\mbox{
\subfigure[a smaller domain]{\includegraphics[scale=0.9]{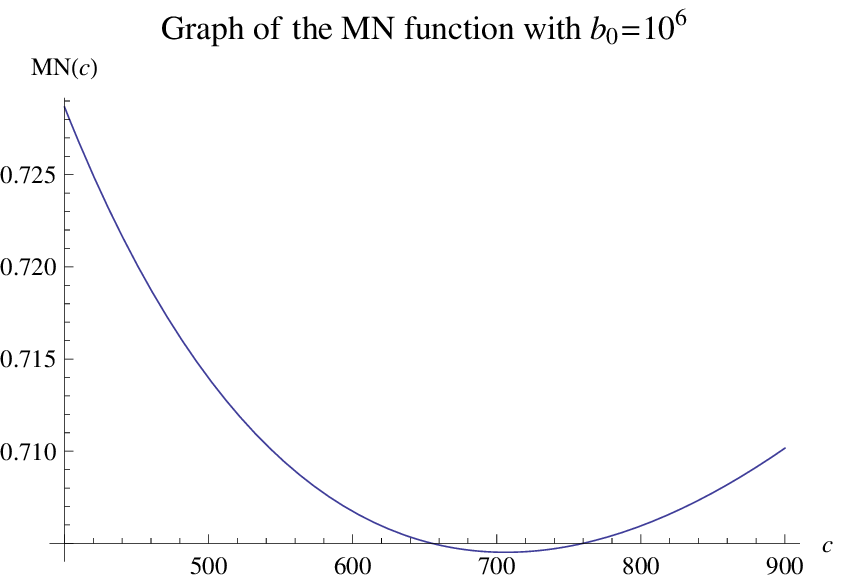}}
\subfigure[a larger domain]{\includegraphics[scale=0.9]{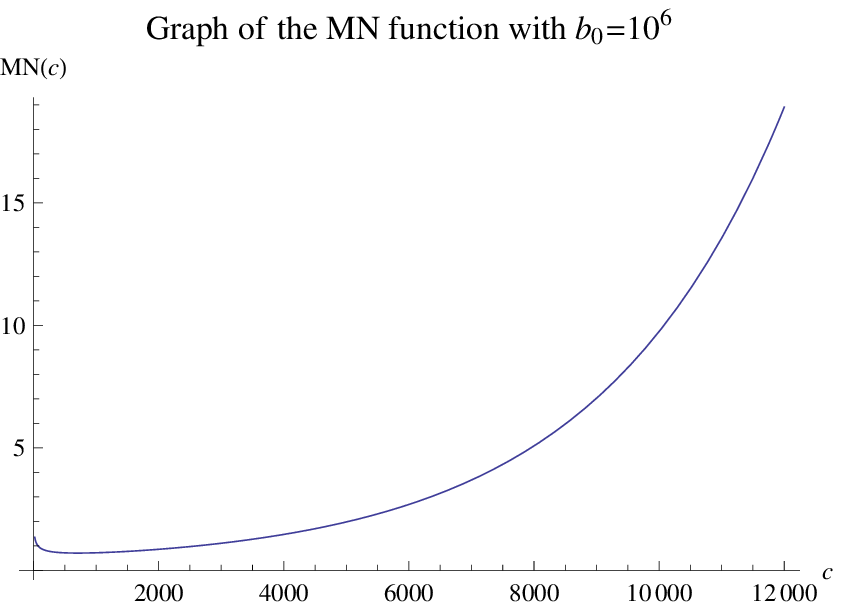}}}
\caption{Here $n=2,\lambda=2,\sigma=1.127$ and $\delta=0.03$.}

\mbox{
\subfigure[a smaller domain]{\includegraphics[scale=0.9]{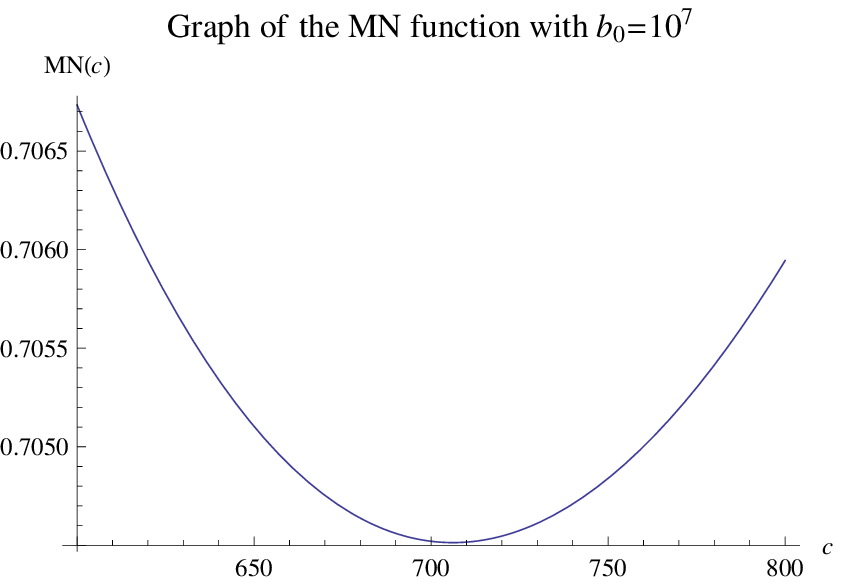}}
\subfigure[a larger domain]{\includegraphics[scale=0.9]{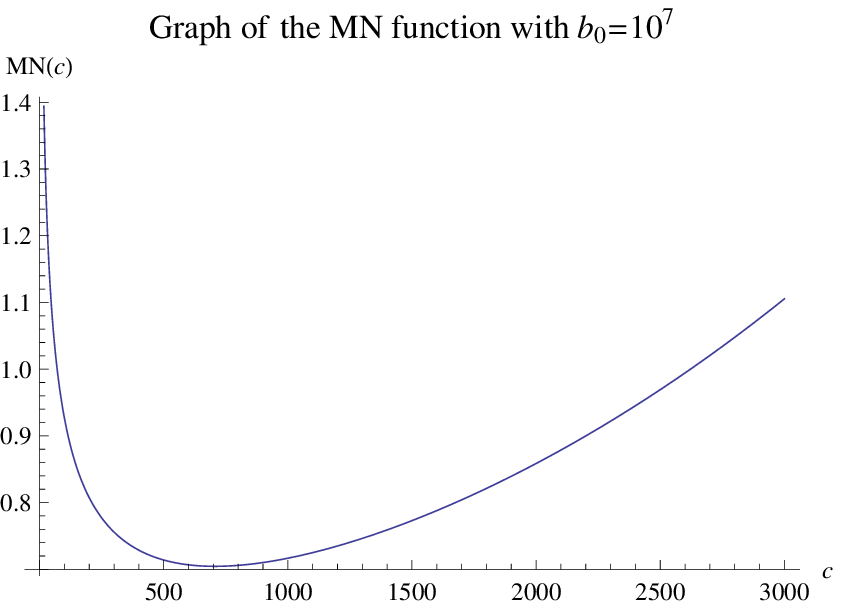}}}
\caption{Here $n=2,\lambda=2,\sigma=1.127$ and $\delta=0.03$.}

\end{figure}
\subsubsection{$f\in E_{\sigma}$} 
For $f\in E_{\sigma}$, there are only two cases.\\
\\
{\bf Case1}. \fbox{$1+\lambda-n\geq 0$} For any $b_{0}>0$ in Theorem2.2 and positive $\delta<\frac{b_{0}}{2(m+1)}$, if $1+\lambda-n\geq 0$, the optimal choice of $c\in [c_{0},\infty)$ is $c^{*}\in [c_{0},c_{1}]$ which minimizes $MN(c)$ in (10) on $[c_{0},c_{1}]$.\\
\\
{\bf Reason}: In this case $MN(c)$ is increasing on $[c_{1},\infty)$. Hence its minimum value happens in $[c_{0},c_{1}]$.\\
\\
{\bf Numerical Examples}:
\begin{figure}[t]
\centering
\mbox{
\subfigure[a smaller domain]{\includegraphics[scale=0.9]{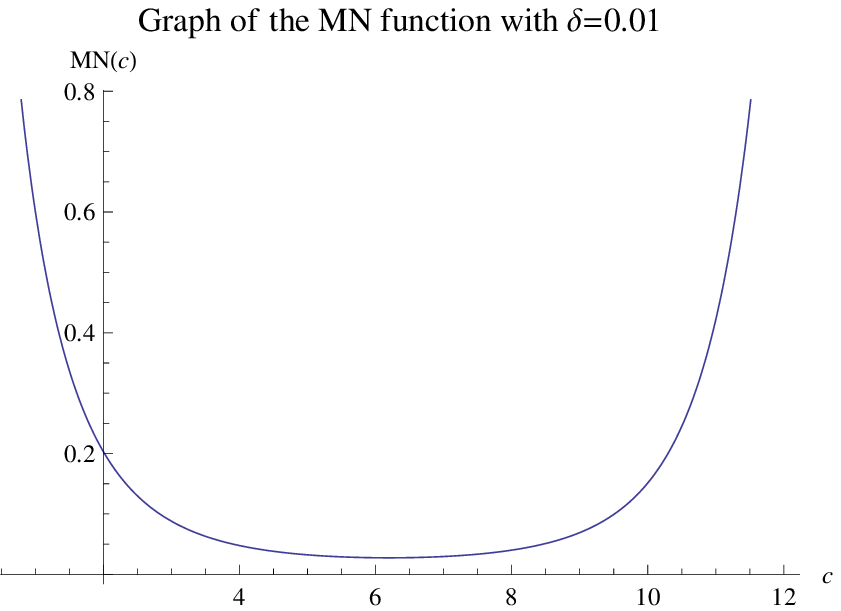}}
\subfigure[a larger domain]{\includegraphics[scale=0.9]{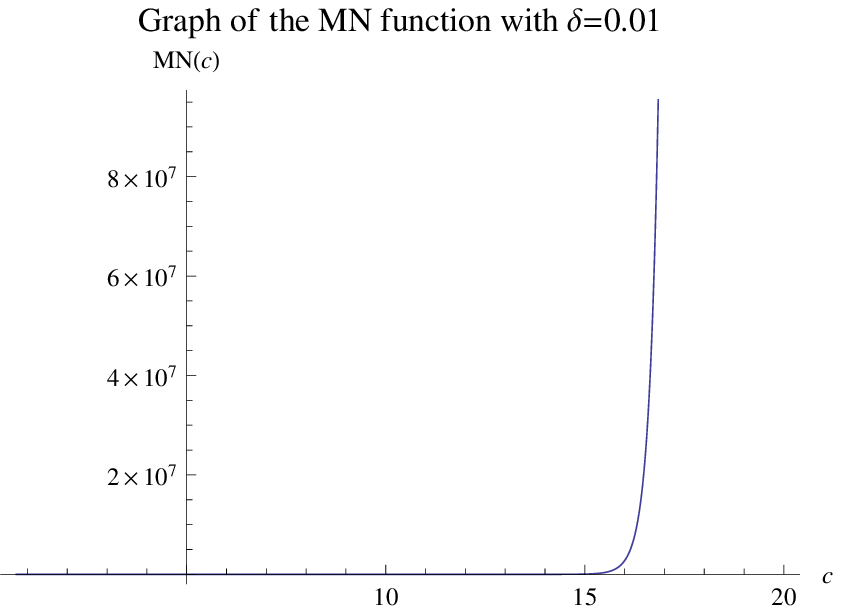}}}
\caption{Here $n=2,\lambda=2,\sigma=1$ and $b_{0}=1$.}

\mbox{
\subfigure[a smaller domain]{\includegraphics[scale=0.9]{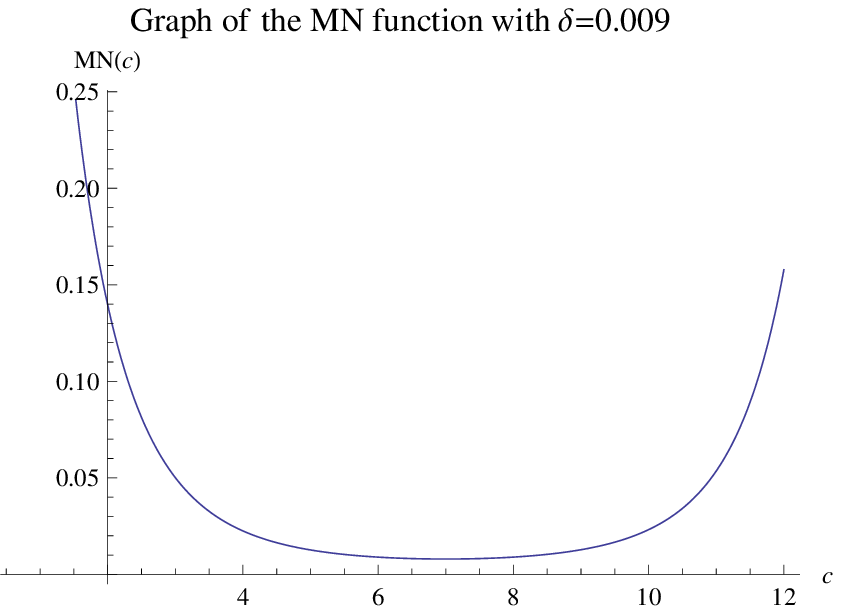}}
\subfigure[a larger domain]{\includegraphics[scale=0.9]{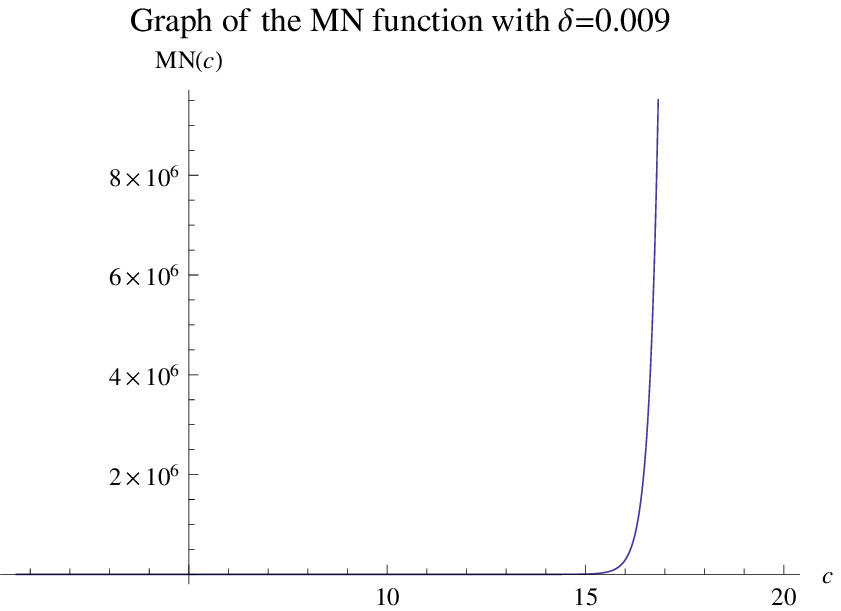}}}
\caption{Here $n=2,\lambda=2,\sigma=1$ and $b_{0}=1$.}

\mbox{
\subfigure[a smaller domain]{\includegraphics[scale=0.9]{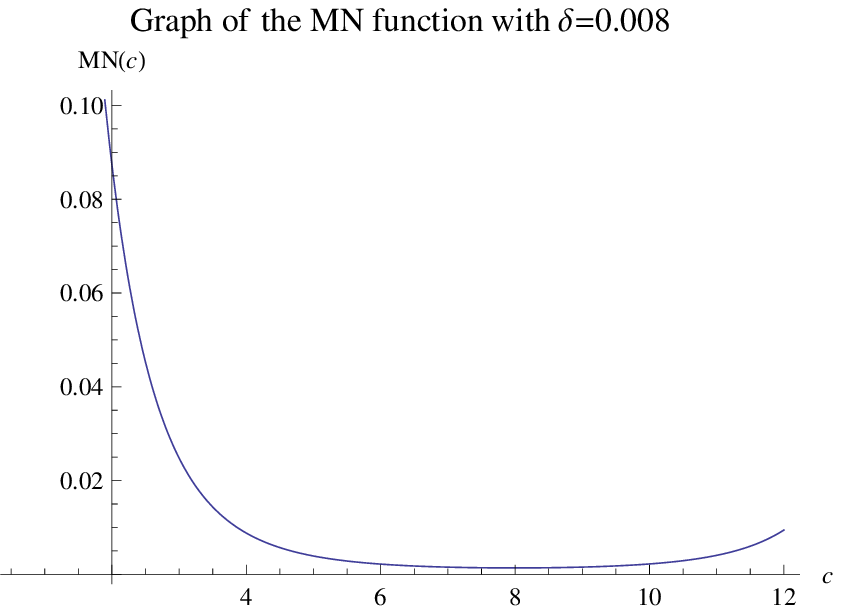}}
\subfigure[a larger domain]{\includegraphics[scale=0.9]{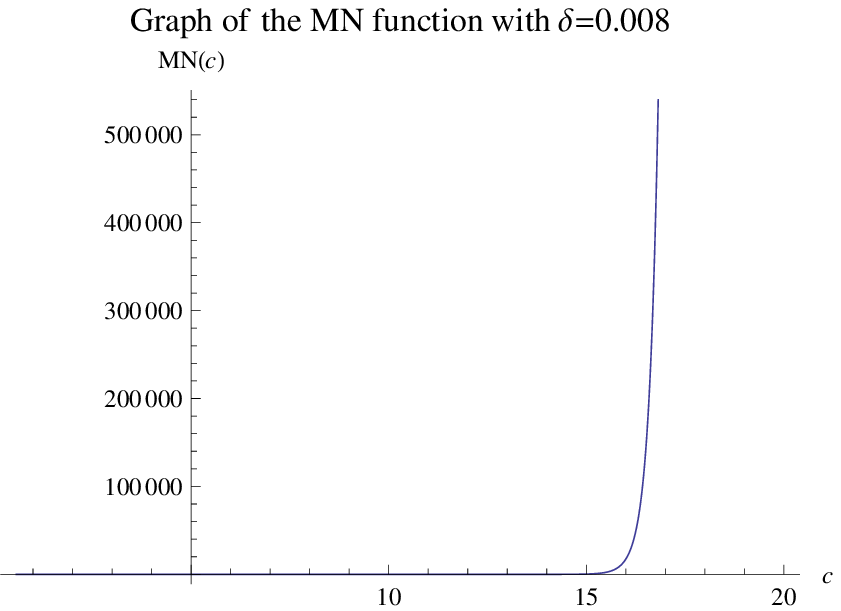}}}
\caption{Here $n=2,\lambda=2,\sigma=1$ and $b_{0}=1$.}

\end{figure}

\clearpage
\begin{figure}[t]
\centering
\mbox{
\subfigure[a smaller domain]{\includegraphics[scale=0.9]{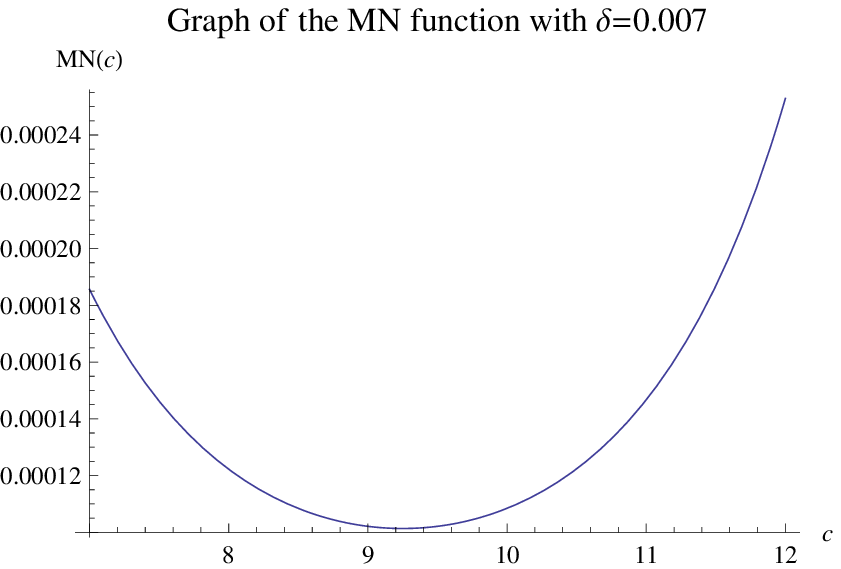}}
\subfigure[a larger domain]{\includegraphics[scale=0.9]{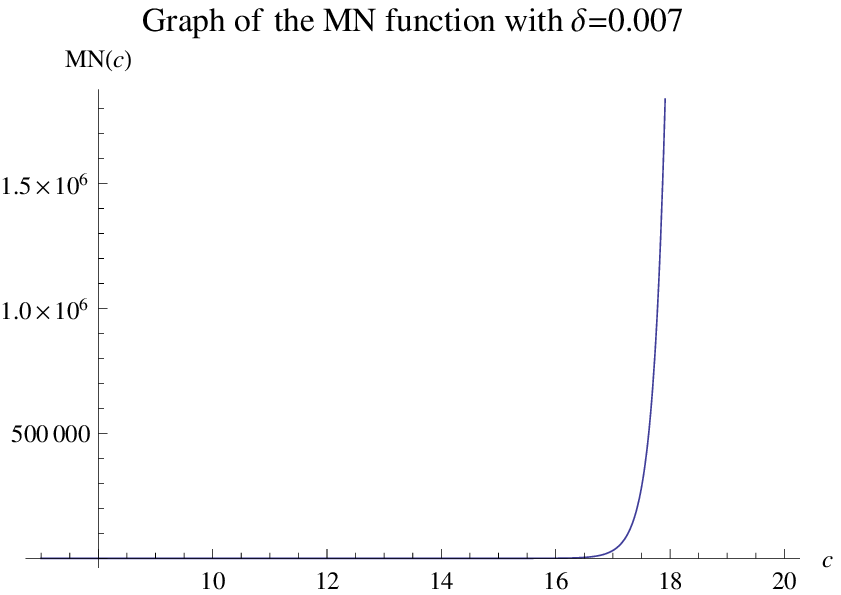}}}
\caption{Here $n=2,\lambda=2,\sigma=1$ and $b_{0}=1$.}

\mbox{
\subfigure[a smaller domain]{\includegraphics[scale=0.9]{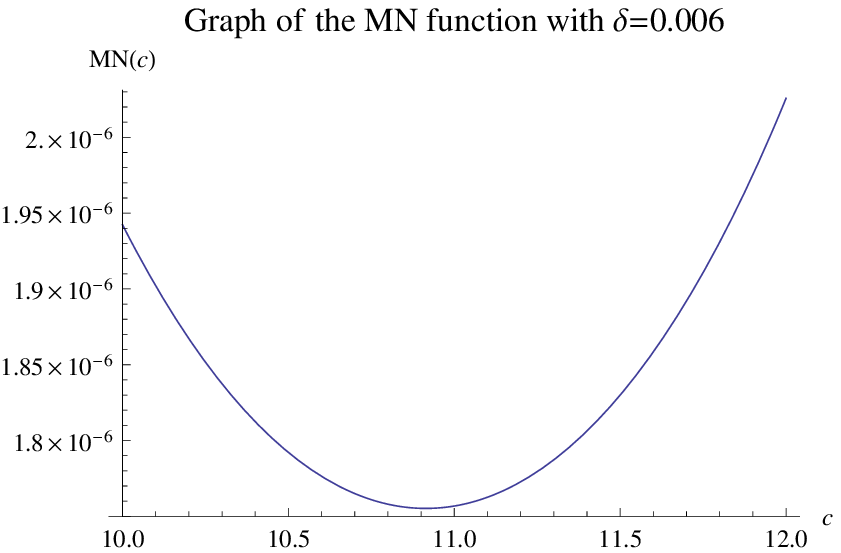}}
\subfigure[a larger domain]{\includegraphics[scale=0.9]{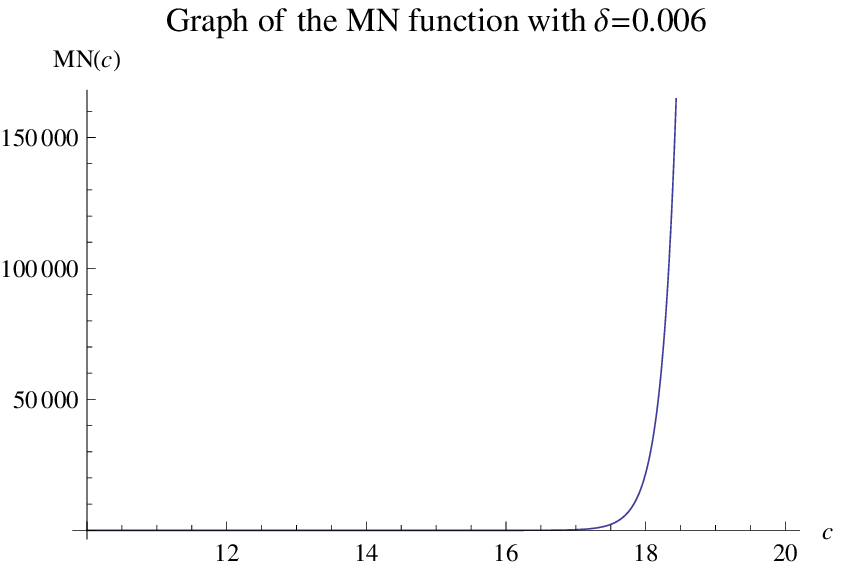}}}
\caption{Here $n=2,\lambda=2,\sigma=1$ and $b_{0}=1$.}

\end{figure}
The following case happens only when $n\geq 4$.\\
\\
{\bf Case2}. \fbox{$1+\lambda-n<0$} For any $b_{0}>0$ in Theorem2.2 and positive $\delta<\frac{b_{0}}{2(m+1)}$, if $1+\lambda-n<0$, the optimal choice of $c\in [c_{0},\infty)$ is either $c^{*}\in[c_{0},c_{1}]$ or $c^{**}\in [c_{1},\infty)$, depending on $MN(c^{*})\leq MN(c^{**})$ or $MN(c^{**})\leq MN(c^{*})$, where $c^{*}$ and $c^{**}$ minimize $MN(c)$ in (10) on $[c_{0},c_{1}]$ and $[c_{1},\infty)$, respectively.\\
\\
{\bf Reason}: In this case $MN(c)$ may not be monotonic on both $[c_{0},c_{1}]$ and $[c_{1},\infty)$.\\
\\
{\bf Remark}: We can apply Matlab or Mathematica to find $c^{*}$ and $c^{**}$.\\
\\
{\bf Numerical Examples}: 
\begin{figure}[t]
\centering
\mbox{
\subfigure[a smaller domain]{\includegraphics[scale=0.9]{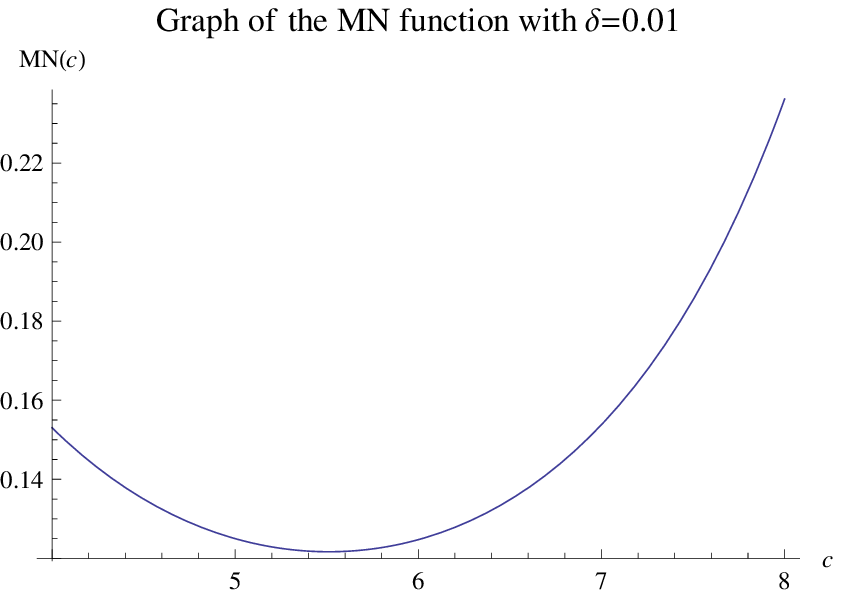}}
\subfigure[a larger domain]{\includegraphics[scale=0.9]{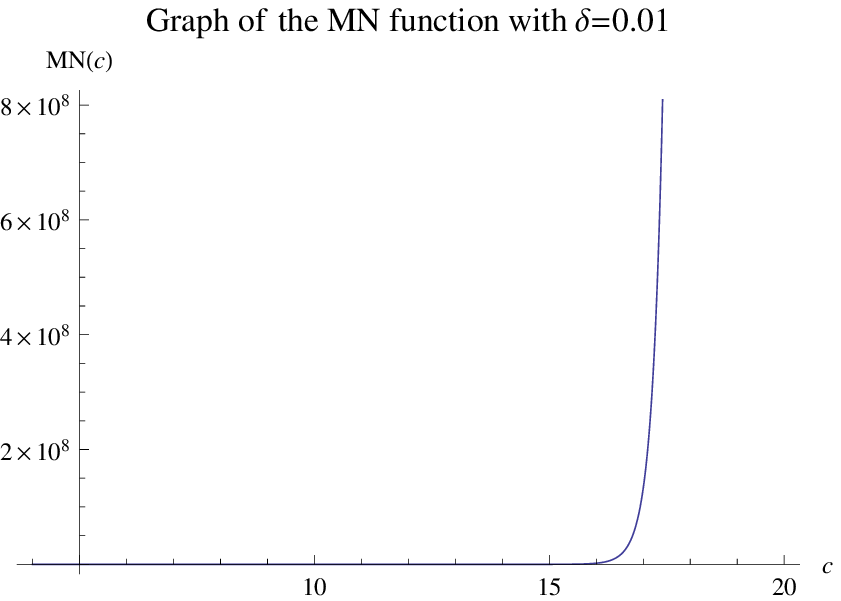}}}
\caption{Here $n=4,\lambda=2,\sigma=1$ and $b_{0}=1$.}

\mbox{
\subfigure[a smaller domain]{\includegraphics[scale=0.9]{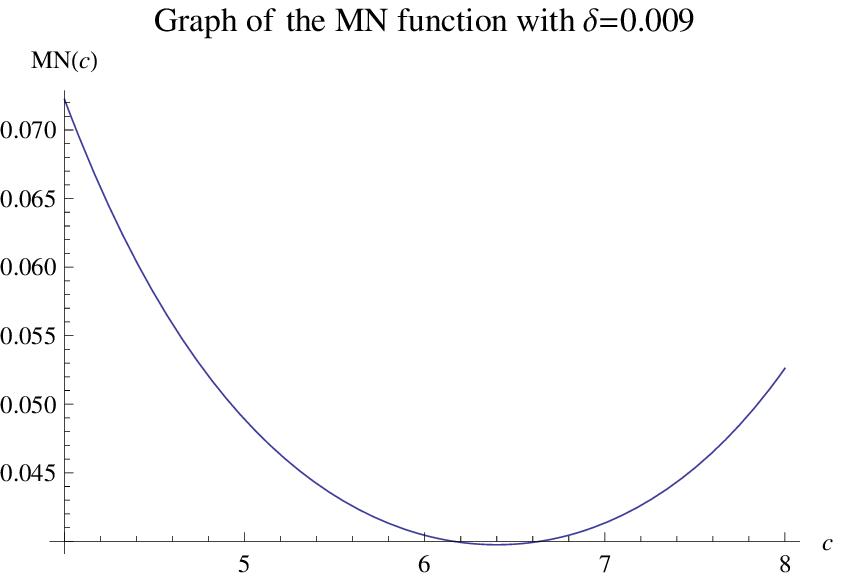}}
\subfigure[a larger domain]{\includegraphics[scale=0.9]{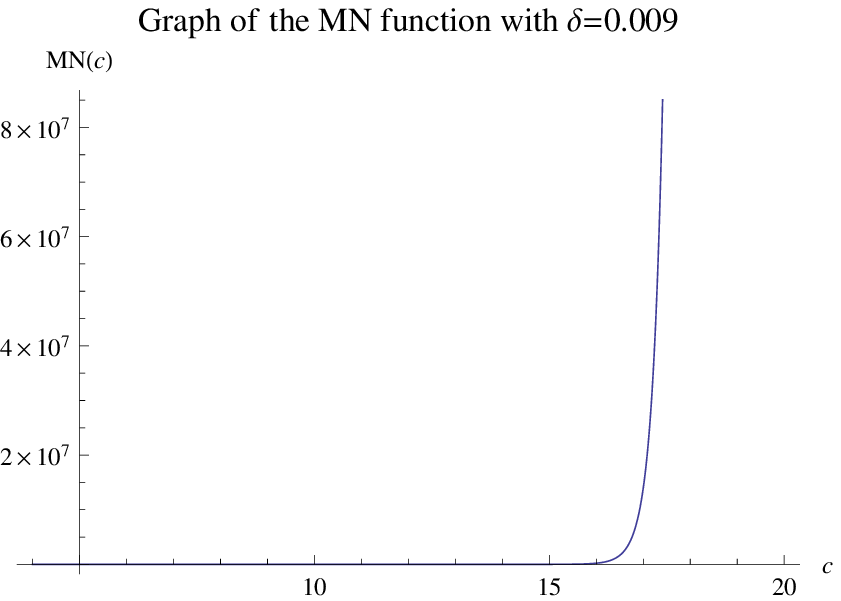}}}
\caption{Here $n=4,\lambda=2,\sigma=1$ and $b_{0}=1$.}

\mbox{
\subfigure[a smaller domain]{\includegraphics[scale=0.9]{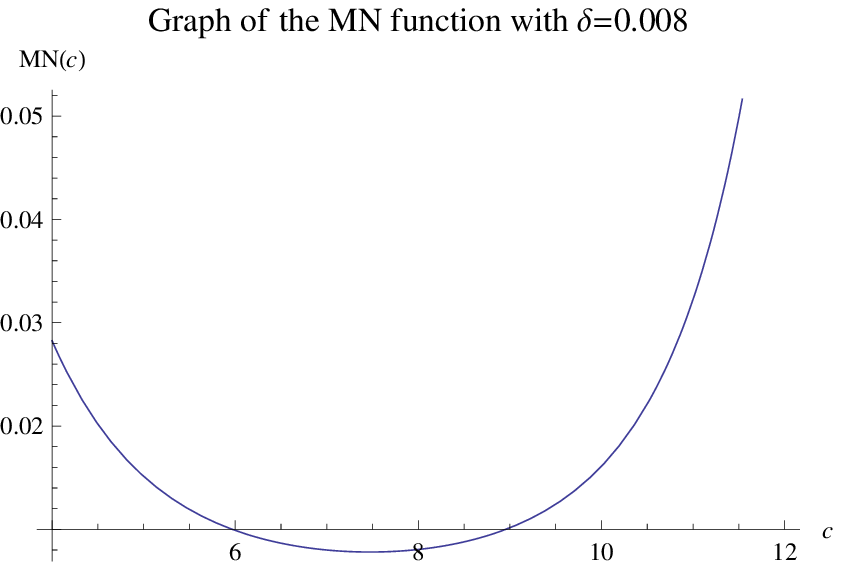}}
\subfigure[a larger domain]{\includegraphics[scale=0.9]{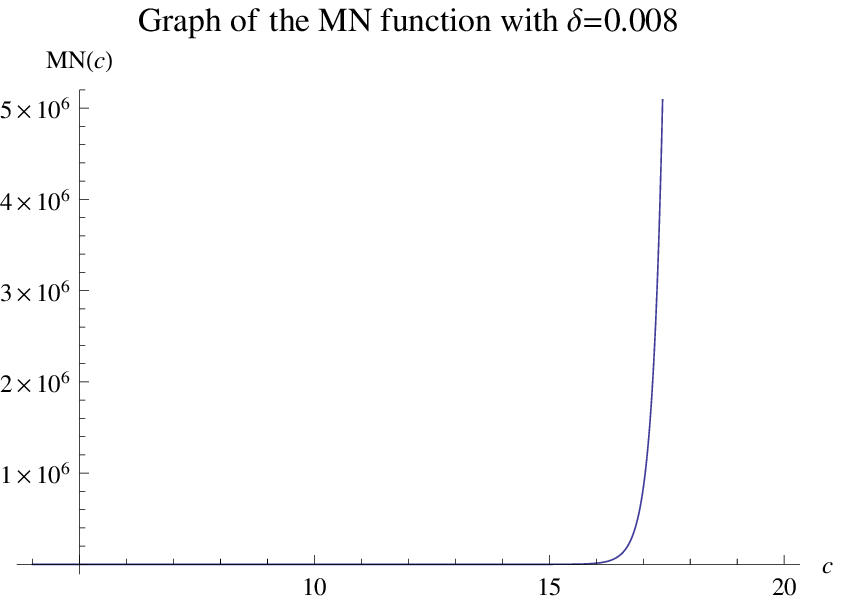}}}
\caption{Here $n=4,\lambda=2,\sigma=1$ and $b_{0}=1$.}

\end{figure}

\clearpage

\begin{figure}[t]
\centering
\mbox{
\subfigure[a smaller domain]{\includegraphics[scale=0.9]{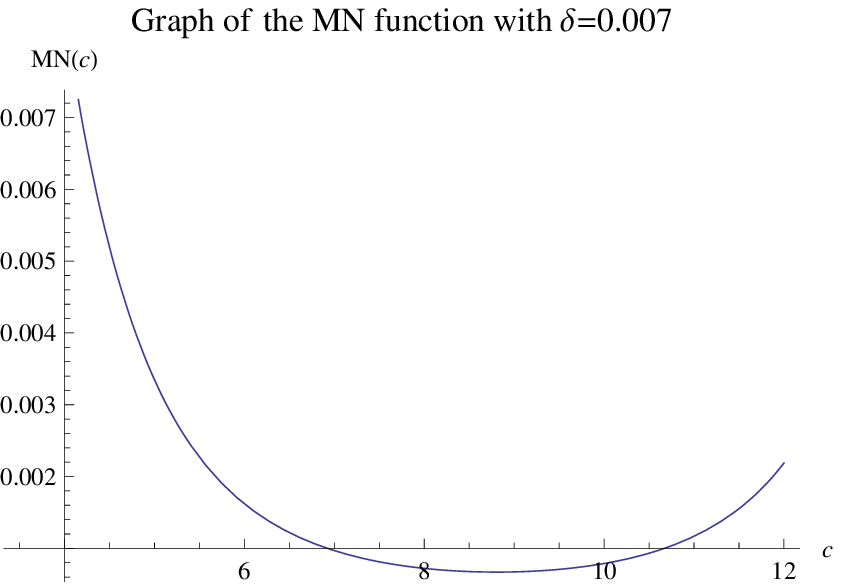}}
\subfigure[a larger domain]{\includegraphics[scale=0.9]{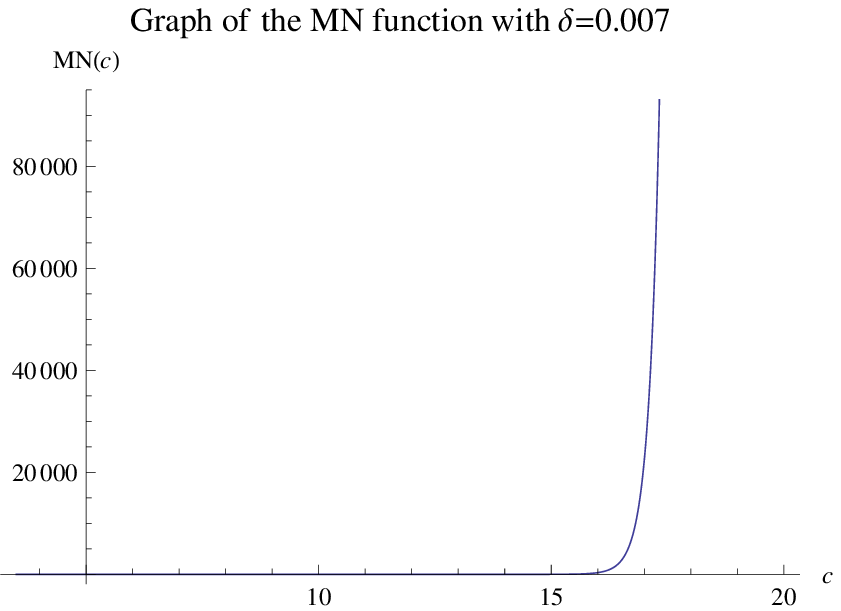}}}
\caption{Here $n=4,\lambda=2,\sigma=1$ and $b_{0}=1$.}

\mbox{
\subfigure[a smaller domain]{\includegraphics[scale=0.9]{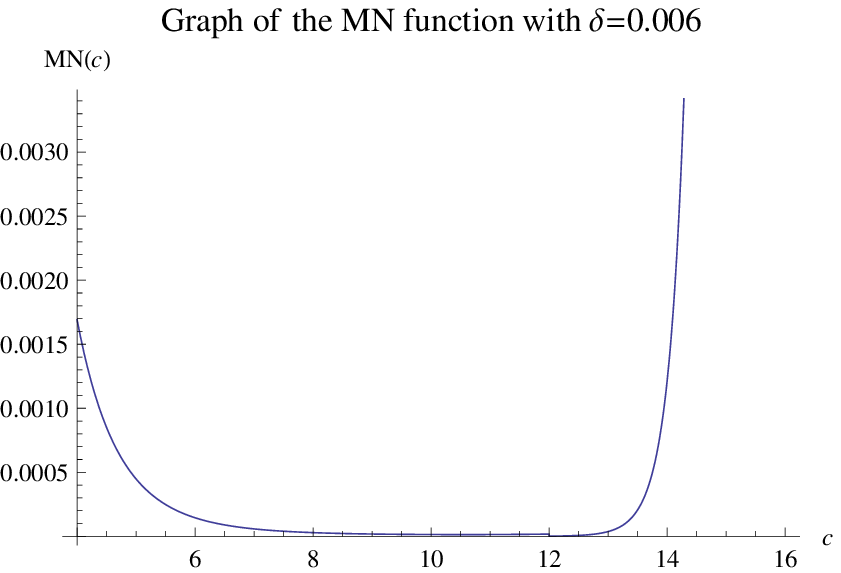}}
\subfigure[a larger domain]{\includegraphics[scale=0.9]{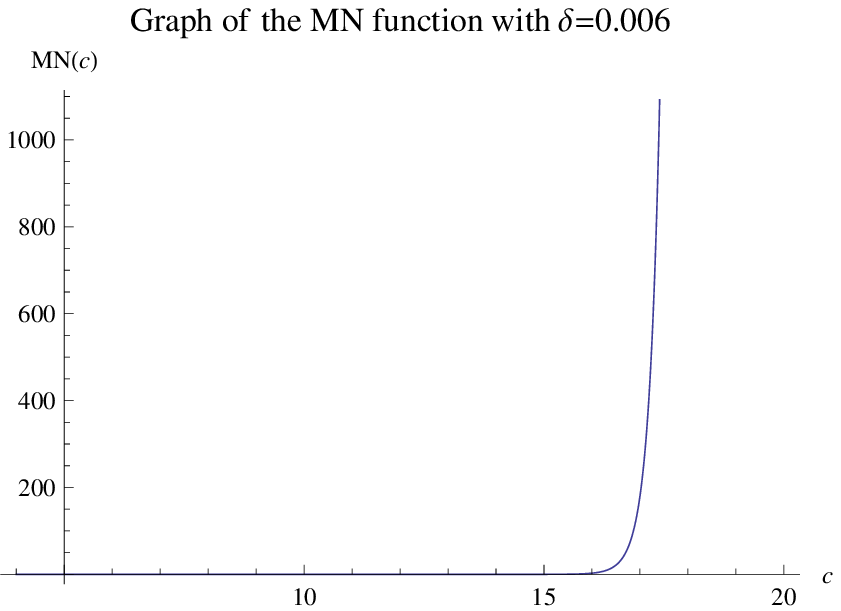}}}
\caption{Here $n=4,\lambda=2,\sigma=1$ and $b_{0}=1$.}
\end{figure}

\subsection{$b_{0}$ not fixed}
For some domains $\Omega$ the number $b_{0}$ in Theorem2.2 can be made arbitrarily large. For example $Q_{0}\subseteq \Omega =R^{n}$ or $Q_{0}\subseteq \Omega =\{ (x_{1},\cdots, x_{n})|\ 0\leq x_{i}<\infty\ for\ i=1,\cdots,n\}$. Such domains $\Omega$ are called dilation-invariant. In this situation one can keep $C=8\rho'=\frac{8\rho}{c}$ in Theorem2.2 by increasing $b_{0}$, while $\delta$ is fixed. It will make the upper bound (4) better because $C$ can become smaller by increasing $c$.  

We never decrease $b_{0}$ because it will make both (4) and the upper bound $\delta_{0}$ of $\delta$ worse.

For any $\delta>0$, once the optimal choice of $c\in [c_{0},\infty)$ is obtained, we let $b_{0}=\frac{c}{12\rho}$. Then $C=8\rho'=8\left(\frac{\rho}{c}\right)=\frac{2}{3b_{0}}$ and $\delta_{0}=\frac{1}{3C(m+1)}=\frac{b_{0}}{2(m+1)}$. Since $c\geq c_{0}=24\rho(m+1)\delta$, we have $\delta \leq \frac{c}{24\rho(m+1)}=\frac{1}{3\left(\frac{8\rho}{c}\right)(m+1)}=\frac{1}{3C(m+1)}=\delta_{0}$, and the requirements of Theorem2.2 are satisfied.

What's noteworthy is that when $b_{0}$ is increased, in order to keep $\delta$ fixed, one has to add more data points to the simplex $Q$ in Theorem2.2.

Now, for any $\delta>0$, the MN functions become
\begin{equation}
 MN(c)=\sqrt{8\rho}\cdot c^{\frac{\lambda-n-1}{4}}e^{c\left[\frac{\sigma}{2}+\frac{\ln{\frac{2}{3}}}{24\rho \delta}\right]},\ c\in [c_{0},\infty )
\end{equation}
for $f\in B_{\sigma}$, and
\begin{equation}
 MN(c)=\sqrt{8\rho}\cdot c^{\frac{\lambda-n-1}{4}}\sup_{\xi\in R^{n}}\left\{ |\xi|^{\frac{1+n+\lambda}{4}}e^{\frac{c|\xi|}{2}-\frac{|\xi|^{2}}{2\sigma}}\right\}\left(\frac{2}{3}\right)^{\frac{c}{24\rho\delta}},\ c\in [c_{0},\infty )
\end{equation}
for $f\in E_{\sigma}$.

Based on (11) and (12), we then have the following criteria of choosing $c$.
\subsubsection{$f\in B_{\sigma}$}
Let $\sigma>0,\ \delta>0$ be fixed, and $k:=\frac{\sigma}{2}+\frac{\ln{\frac{2}{3}}}{24\rho\delta}$. Let $f\in B_{\sigma}$ be the approximated function in Theorem2.2.\\
\\
{\bf Case1}. \fbox{$\lambda-n-1\geq 0$} and \fbox{$k>0$} If $\lambda-n-1\geq 0$ and $k>0$, the optimal choice of $c\in [c_{0},\infty )$ is $c=c_{0}:=24\rho(m+1)\delta$.\\
\\
{\bf Reason}: In this case $MN(c)$ in (11) is increasing on $[c_{0}, \infty )$.\\
\\
{\bf Case2}. \fbox{$\lambda-n-1\leq 0$} and \fbox{$k<0$} If $\lambda-n-1\leq 0$ and $k<0$, the larger $c$ is, the better it is.\\
\\
{\bf Reason}: In this case $MN(c)$ in (11) is decreasing and $MN(c)\rightarrow 0$ as $c\rightarrow \infty$.\\
\\
{\bf Case3}. If $\lambda-n-1$ and $k$ are of opposite signs, $MN(c)$ in (11) is not monotonic and the optimal $c\in [c_{0},\infty )$ is the number $c^{*}$ minimizing $MN(c)$, which can be found by Matlab or Mathematica.
{\bf Numerical Examples}:
\begin{figure}[h]
\centering
\mbox{
\subfigure[$\sigma=3.37898$]{\includegraphics[scale=0.9]{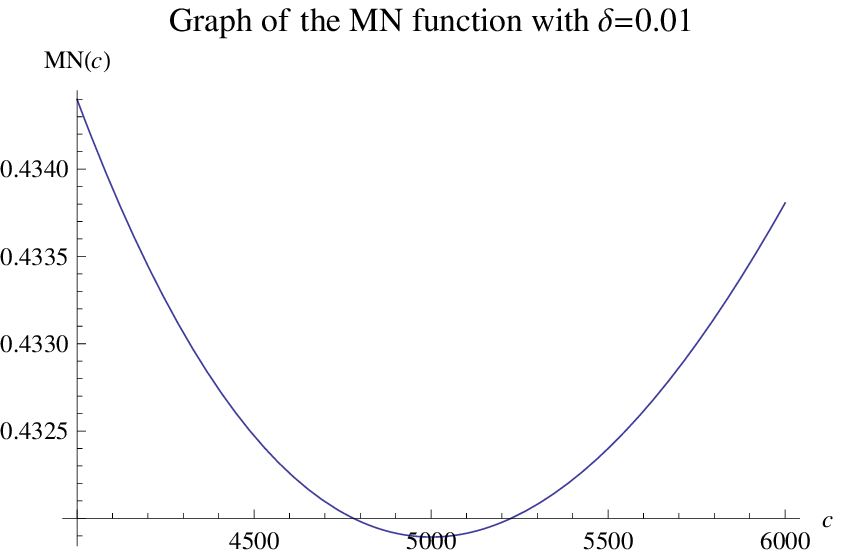}}
\subfigure[$\sigma=3.37908$]{\includegraphics[scale=0.9]{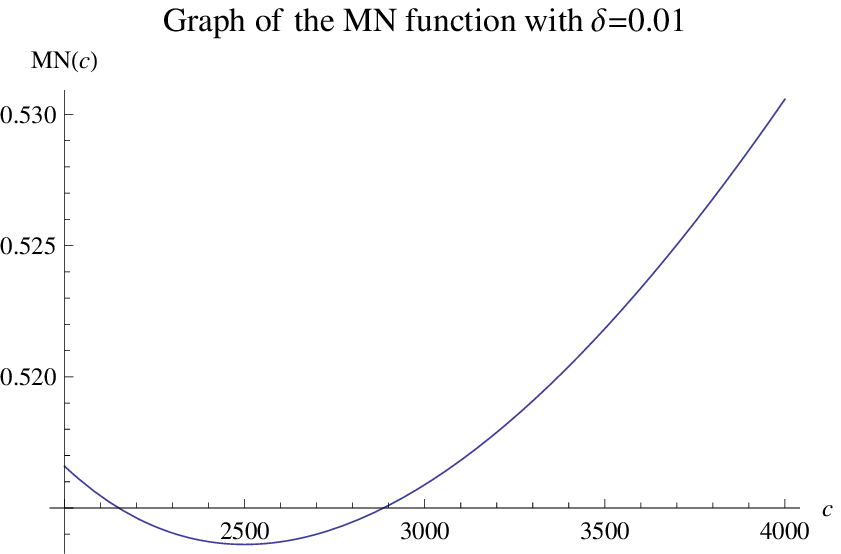}}}
\caption{Here $n=2$ and $\lambda=2$.}
\end{figure}

\clearpage
\begin{figure}[t]
\centering
\mbox{
\subfigure[$\sigma=4.25431$]{\includegraphics[scale=0.9]{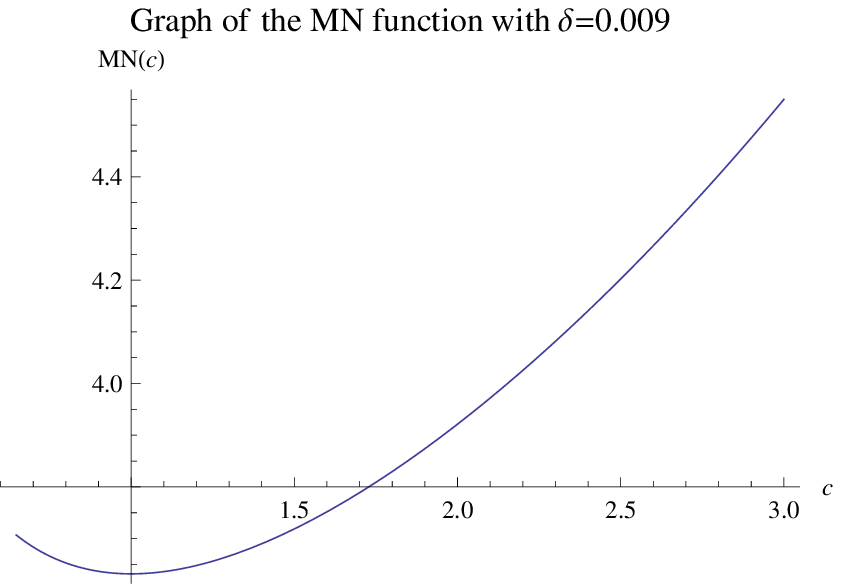}}
\subfigure[$\sigma=3..95431$]{\includegraphics[scale=0.9]{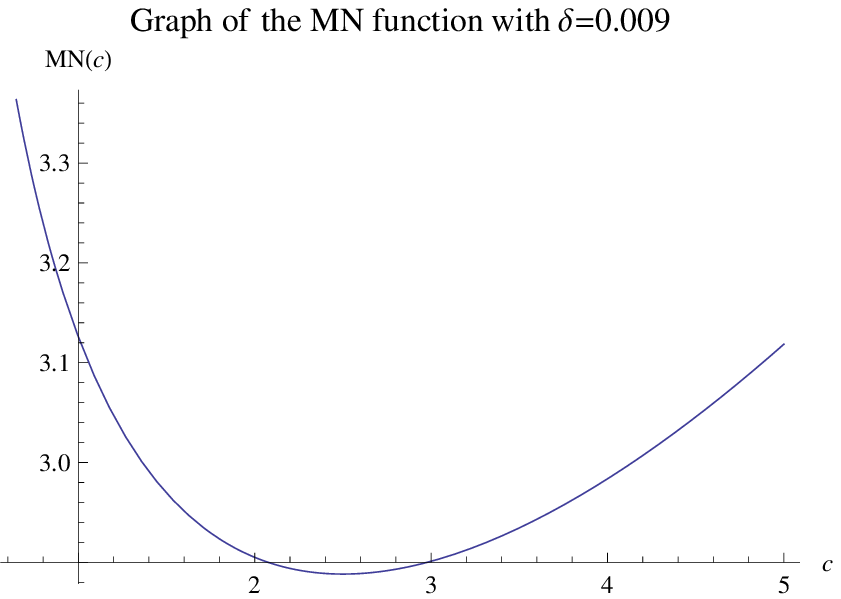}}}
\caption{Here $n=2$ and $\lambda=2$.}

\mbox{
\subfigure[$\sigma=4.72359$]{\includegraphics[scale=0.9]{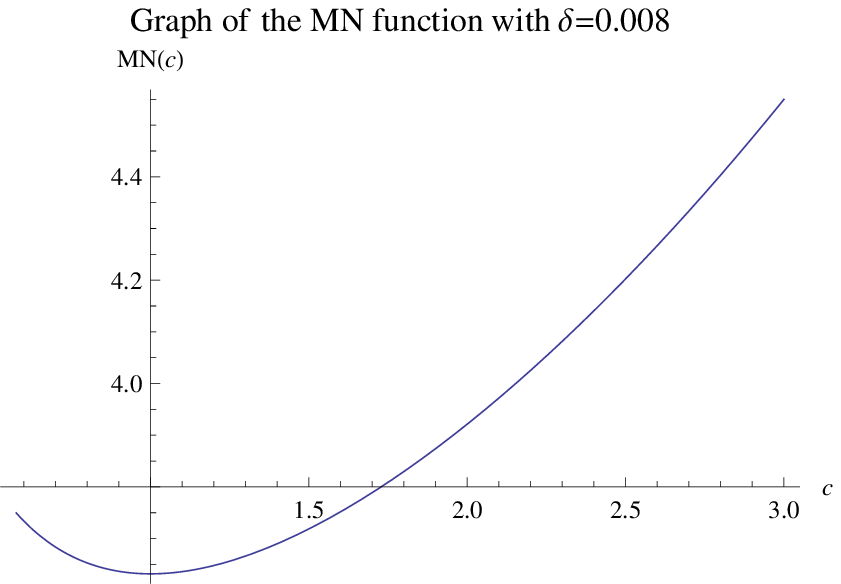}}
\subfigure[$\sigma=4.42359$]{\includegraphics[scale=0.9]{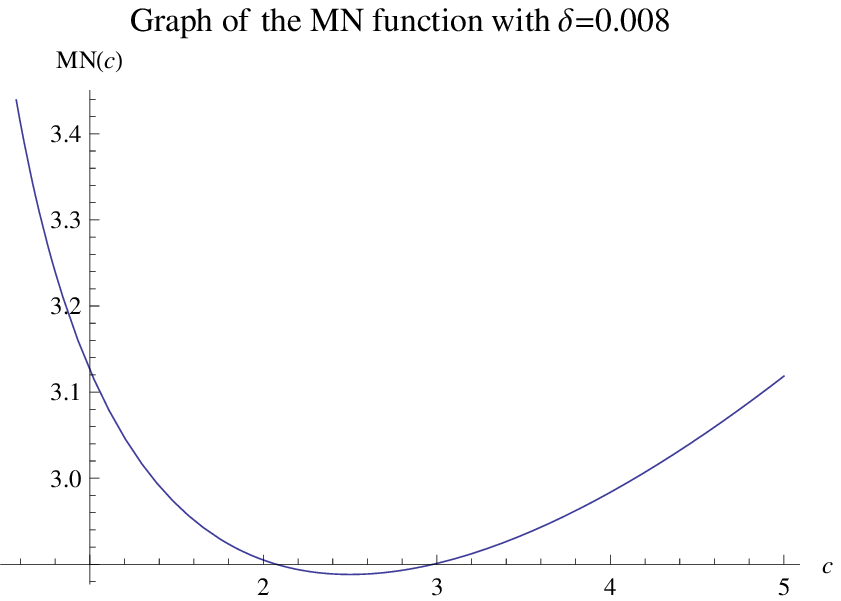}}}
\caption{Here $n=2$ and $\lambda=2$.}

\mbox{
\subfigure[$\sigma=5.32697$]{\includegraphics[scale=0.9]{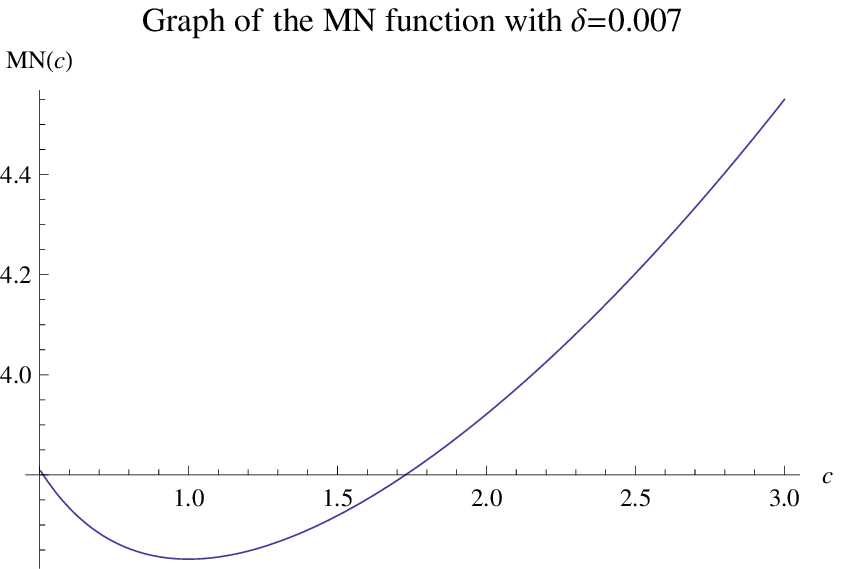}}
\subfigure[$\sigma=5.02697$]{\includegraphics[scale=0.9]{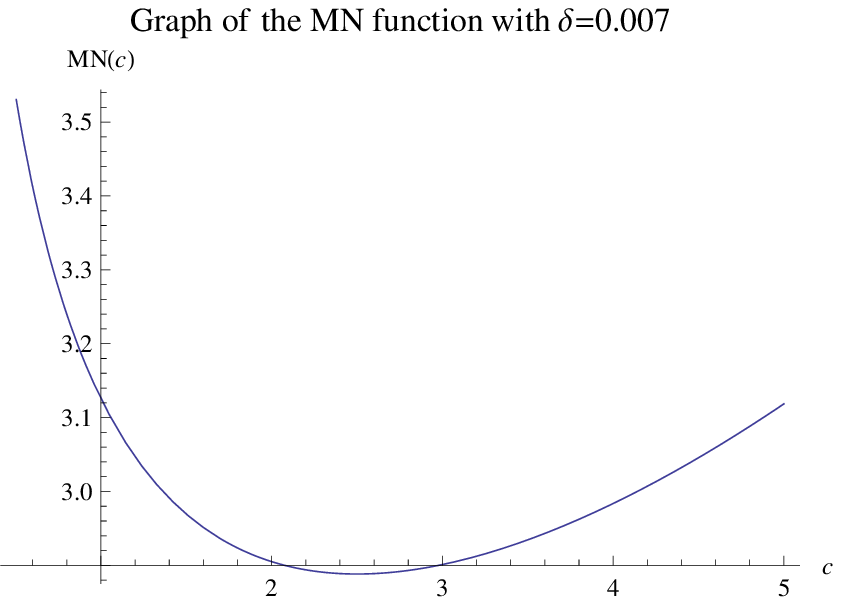}}}
\caption{Here $n=2$ and $\lambda=2$.}
\end{figure}

\clearpage
\begin{figure}[t]
\centering
\mbox{
\subfigure[$\sigma=6.13146$]{\includegraphics[scale=0.9]{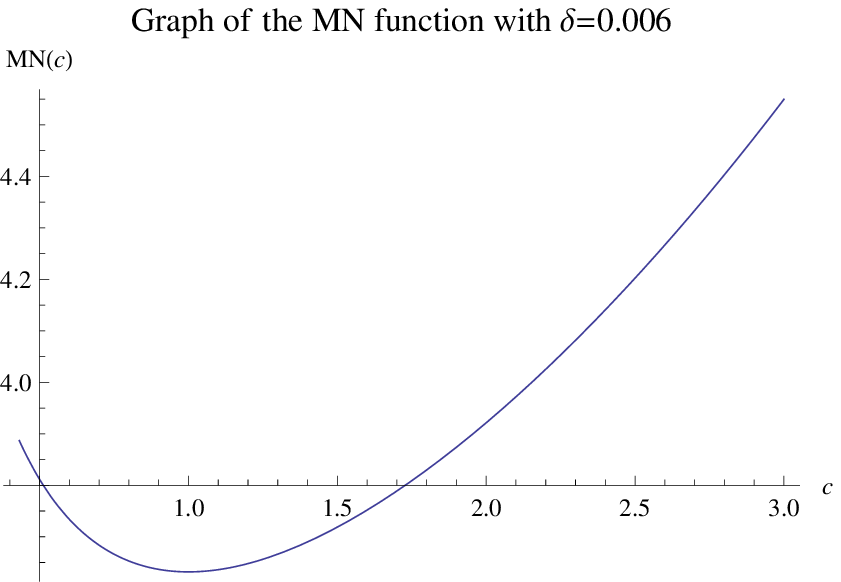}}
\subfigure[$\sigma=5.83146$]{\includegraphics[scale=0.9]{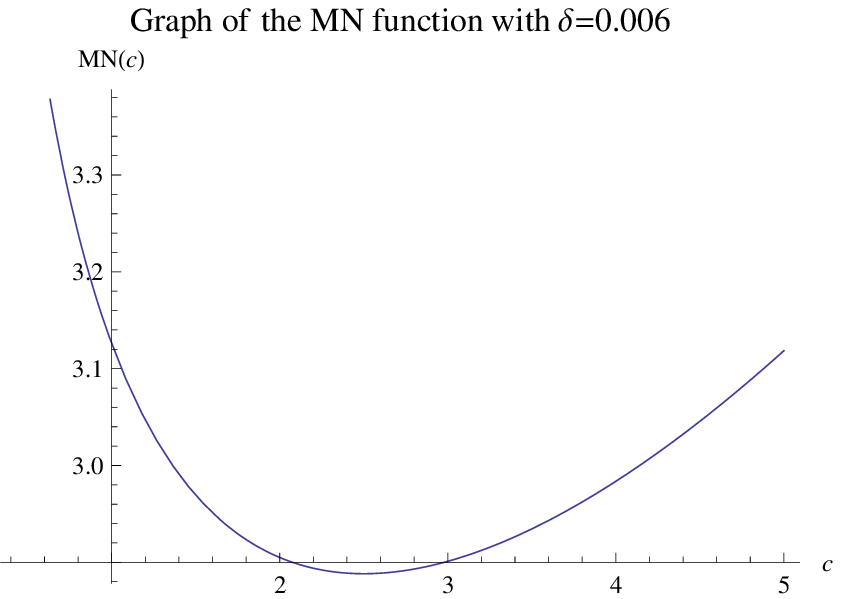}}}
\caption{Here $n=2$ and $\lambda=2$.}
\end{figure}
Note that:(a)in figure27 the two $\sigma$'s differ only by 0.0001. In Figure28-31 the difference is only 0.3;(b)in the preceding five examples the optimal choice of $c$ is very sensitive to $\sigma$, but not to $\delta$.
\subsubsection{$f\in E_{\sigma}$}
For $f\in E_{\sigma}$, $MN(c)$ in (12) is not monotonic. Hence one should use Matlab or Mathematica to find $c^{*}\in [c_{0},\infty)$ which minimizes $MN(c)$.\\
\\
{\bf Numerical Examples}:
\begin{figure}[h]
\centering
\mbox{
\subfigure[a smaller domain]{\includegraphics[scale=0.9]{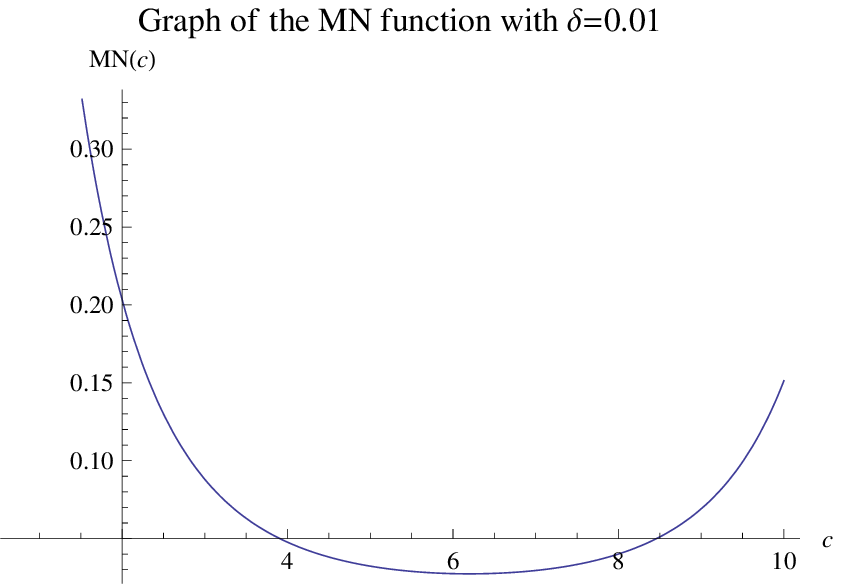}}
\subfigure[a larger domain]{\includegraphics[scale=0.9]{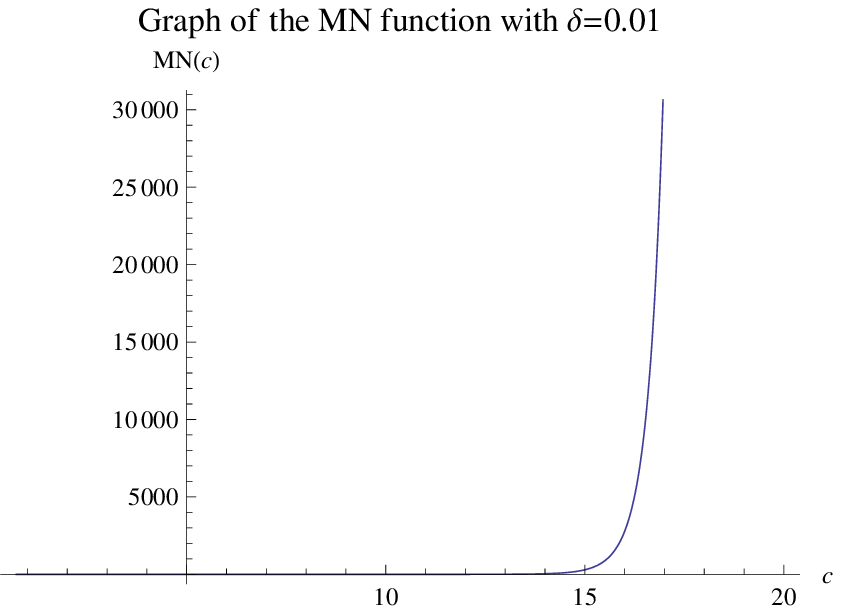}}}
\caption{Here $n=2,\lambda=2$ and $\sigma=1$.}
\end{figure}

\clearpage

\begin{figure}[t]
\centering
\mbox{
\subfigure[a smaller domain]{\includegraphics[scale=0.9]{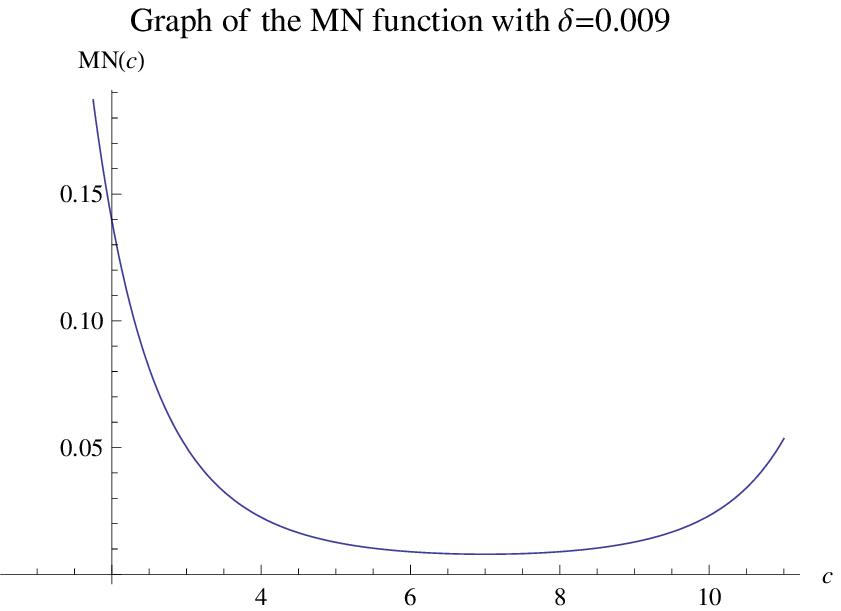}}
\subfigure[a larger domain]{\includegraphics[scale=0.9]{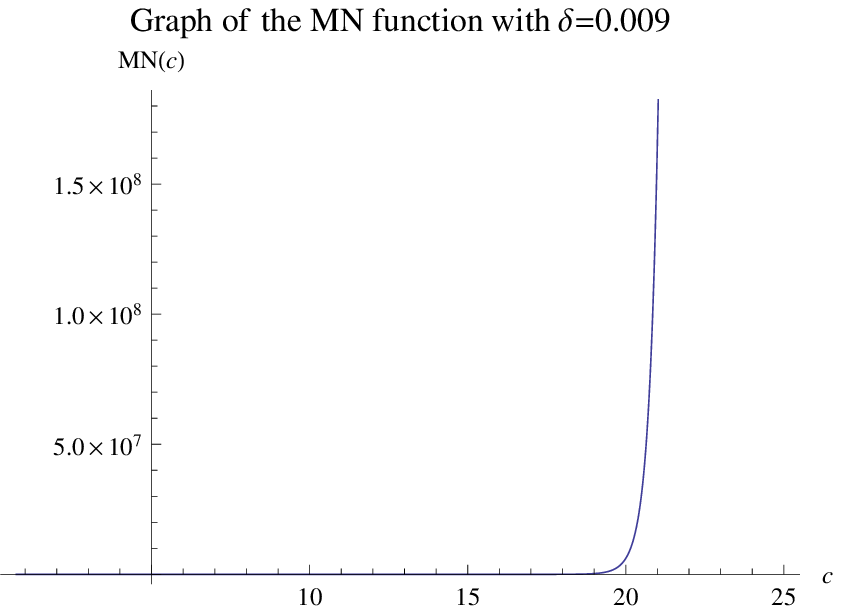}}}
\caption{Here $n=2,\lambda=2$ and $\sigma=1$.}

\mbox{
\subfigure[a smaller domain]{\includegraphics[scale=0.9]{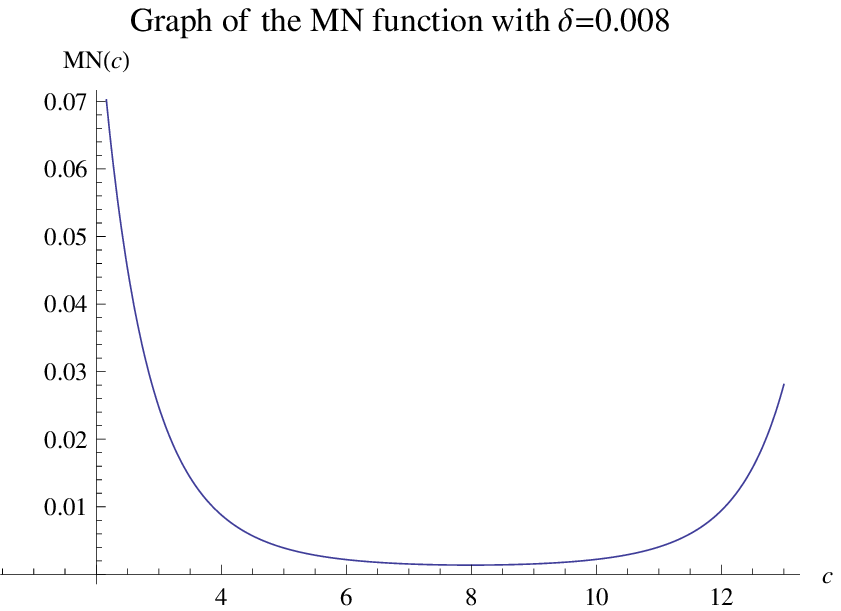}}
\subfigure[a larger domain]{\includegraphics[scale=0.9]{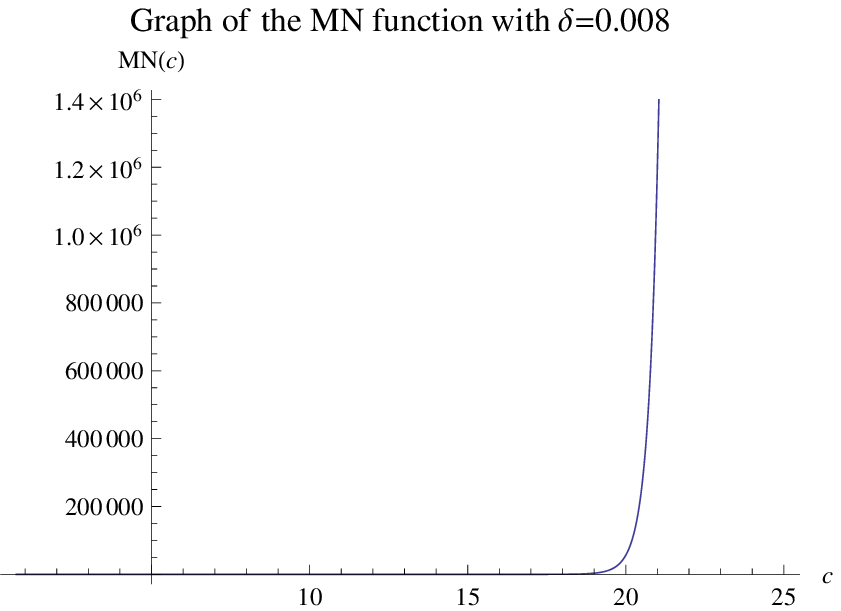}}}
\caption{Here $n=2,\lambda=2$ and $\sigma=1$.}

\mbox{
\subfigure[a smaller domain]{\includegraphics[scale=0.9]{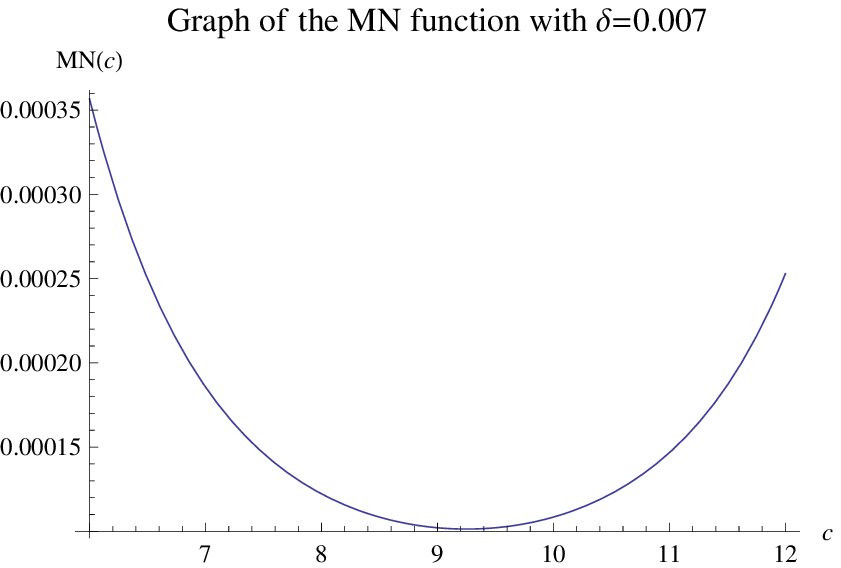}}
\subfigure[a larger domain]{\includegraphics[scale=0.9]{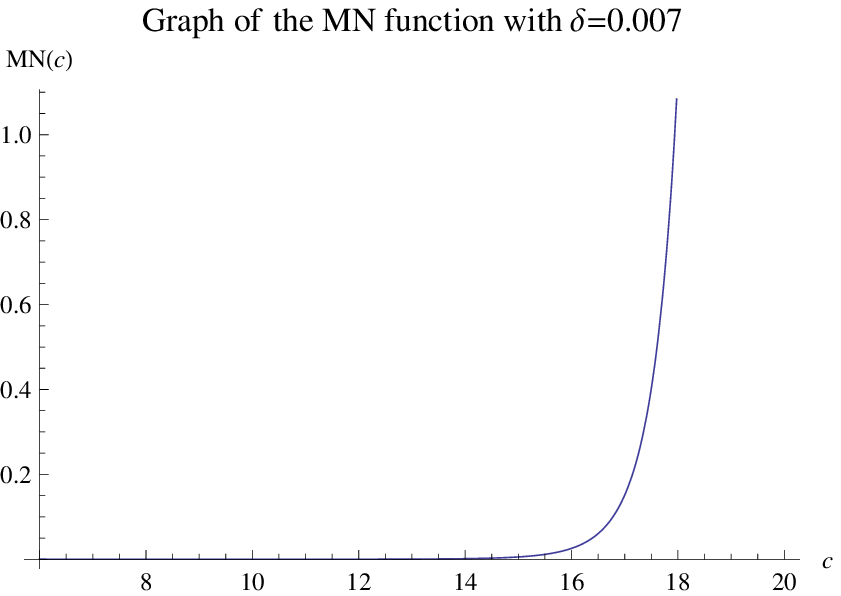}}}
\caption{Here $n=2,\lambda=2$ and $\sigma=1$.}
\end{figure}

\clearpage

\begin{figure}[t]
\centering
\mbox{
\subfigure[a smaller domain]{\includegraphics[scale=0.9]{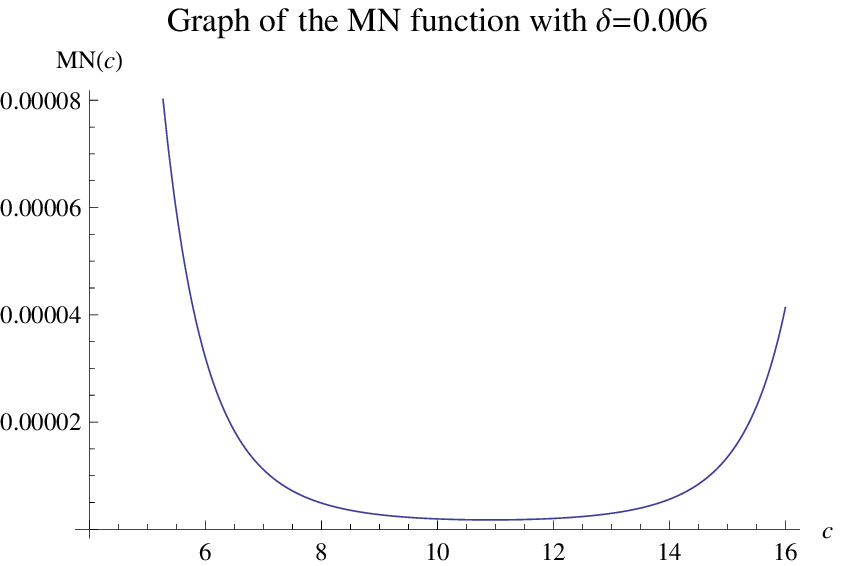}}
\subfigure[a larger domain]{\includegraphics[scale=0.9]{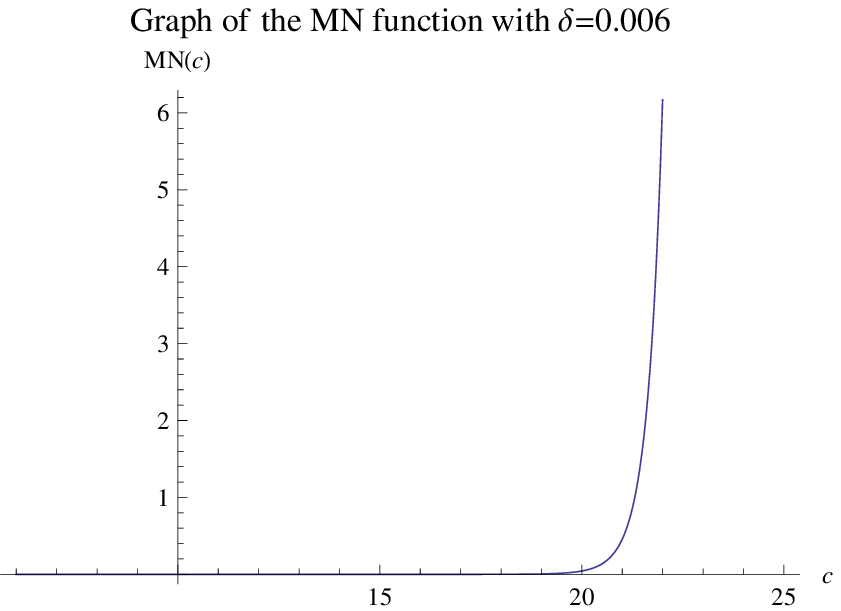}}}
\caption{Here $n=2,\lambda=2$ and $\sigma=1$.}
\end{figure}


\end{document}